\documentclass[12pt]{amsart}
\usepackage{amssymb}

\newtheorem{theorem}{Theorem}[section]
\newtheorem{corollary}[theorem]{Corollary}
\newtheorem{proposition}[theorem]{Proposition}

\newtheorem{lemma}[theorem]{Lemma}
\theoremstyle{definition}    
\newtheorem{definition}[theorem]{Definition}
\theoremstyle{remark}

\newtheorem{remark}[theorem]{Remark}

\newtheorem{example}[theorem]{Example}
\newcommand{\beq}{\begin{eqnarray*}}
\newcommand{\eeq}{\end{eqnarray*}}

\newcommand{\D}{\mathfrak{D}}
\newcommand\A{\mathcal{A}}

\renewcommand{\L}{\mathcal{L}}

\renewcommand{\O}{\mathcal{O}}

\newcommand{\Co}{\mathcal{C}}
\newcommand{\ca}{\mathcal}

\newcommand{\R}{\mathbb{R}}
\newcommand{\C}{\mathbb{C}}

\newcommand{\Z}{\mathbb{Z}}

\newcommand\lie[1]{\mathfrak{#1}}

\newcommand{\h}{\lie{h}}
\newcommand{\g}{\lie{g}}

\renewcommand{\t}{\lie{t}}

\newcommand{\on}{\operatorname}

\newcommand{\Ad}{ \on{Ad} }
\newcommand{\ad}{\on{ad}}
\newcommand{\Hom}{ \on{Hom}}

\renewcommand{\ker}{ \on{ker}}
\newcommand{\im}{ \on{im}}

\newcommand{\Mult}{  \on{Mult}}

\newcommand{\diag}{  \on{diag}}
\newcommand\qu{/\kern-.7ex/} % Categorical quotients
\newcommand{\fus}{\circledast} % Fusion product
\newcommand{\hra}{\hookrightarrow}

\renewcommand{\d}{{\mbox{d}}}

\newcommand\Sig{\Sigma}

\newcommand\Om{\Omega}
\newcommand\om{\omega}

\newcommand{\f}{\frac}

\newcommand{\p}{\partial}
\renewcommand{\l}{\langle}
\renewcommand{\r}{\rangle}
\newcommand\hh{{\f{1}{2}}}

\newcommand{\Id}{\on{Id}}
\newcommand{\wh}{\widehat}
\newcommand{\wt}{\widetilde}
\newcommand{\mf}{\mathfrak}

\begin{document}

\title{Quasi-Poisson manifolds}

\author[A. Alekseev]{A. Alekseev}
\address{Institute for Theoretical Physics \\ Uppsala University \\
Box 803 \\ \mbox{S-75108} Uppsala \\ Sweden}
\email{alekseev@teorfys.uu.se}

\author[Y. Kosmann-Schwarzbach]{Y. Kosmann-Schwarzbach}
\address{Centre de Math\'ematiques (U.M.R. du C.N.R.S. 7640) \\
Ecole Polytechnique, F-91128 Palaiseau, France}
\email{yks@math.polytechnique.fr}

\author[E. Meinrenken]{E. Meinrenken}
\address{University of Toronto, Department of Mathematics,
100 St George Street, Toronto, Ontario M5S3G3, Canada }
\email{mein@math.toronto.edu}

\begin{abstract}
A quasi-Poisson manifold is a $G$-manifold equipped with an invariant 
bivector field whose Schouten bracket is the trivector field 
generated by the invariant element in $\wedge^3\g$ associated to an
invariant inner product. We 
introduce the concept of the fusion of such manifolds, and we relate 
the quasi-Poisson manifolds to the previously introduced
quasi-Hamiltonian manifolds with group-valued moment maps. 
\end{abstract}

\maketitle

\section{Introduction}
This paper is a sequel to \cite{al:ma}, in which the notion 
of a {\em quasi-Poisson manifold} was introduced. While the purpose 
of \cite{al:ma} was to obtain a unified picture of various 
notions of ``generalized moment maps'' and their properties, 
the current article is devoted to a closer examination of a particular
type of quasi-Poisson manifolds, defined as follows.

Let $G$ be a Lie group, who
se Lie algebra $\g$ is equipped
with an invariant inner product. 
 From the Lie bracket 
and the invariant inner 
product, one obtains an invariant element $\phi\in\wedge^3\g$, 
the Cartan $3$-tensor of the Lie algebra with an invariant inner product. 
The quasi-Poisson manifolds studied in 
the present paper are $G$-manifolds $M$ together with an invariant 
bivector field $P$, such that the Schouten 
bracket $[P,P]$ equals the trivector field $\phi_M$ generated 
by $\phi$. Of particular interest are the quasi-Poisson manifolds 
admitting a {\em moment map}, $\Phi:\,M\to G$ (see Definition 
\ref{def:quasiPoisson} below). The triple $(M,P,\Phi)$ will then be called
a {\it Hamiltonian quasi-Poisson manifold}. 

The basic example of a Hamiltonian quasi-Poisson $G$-manifold
is $M=G$, where $G$ acts by conjugation, and the moment map is the identity 
map. As we will explain in this paper, there are many parallels 
between this example and the usual linear Poisson structure on the 
dual of a Lie algebra. In particular, in analogy with the
coadjoint orbits,
all the conjugacy classes in $G$ 
are quasi-Poisson submanifolds with a ``non-degenerate'' 
bivector field. Other examples are obtained using the methods
of ``fusion'' and ``exponentiation'' introduced in this paper.   
As an application, we will show how to construct the usual Poisson structure 
on the representation variety, $\Hom(\pi_1(\Sig),G)/G$, for 
any oriented surface with boundary, $\Sig$, by reduction of a 
quasi-Poisson structure on a fusion product of several copies of $G$. 

One of the main results of this paper states that every Hamiltonian 
quasi-Poisson manifold has a generalized foliation, the leaves of which are 
non-degenerate Hamiltonian quasi-Poisson manifolds. Furthermore, 
every non-degenerate Hamiltonian quasi-Poisson manifold $(M,P,\Phi)$
carries an invariant 2-form $\om$ such that $(M,\om,\Phi)$ satisfies the 
axioms of a {\it quasi-Hamiltonian}
$G$-{\it manifold} with group-valued moment map, 
as introduced in \cite{al:mom}. 
Conversely, every such manifold carries a non-degenerate
quasi-Poisson structure with moment map $\Phi$. 

In a recent paper, T. Treloar \cite{tr:sy} has applied the concept 
of a Hamiltonian quasi-Poisson manifold developed here 
to study the 
symplectic geometry of spaces of polygons 
on the 3-sphere. 

The organization of the paper is as follows.
In Section \ref{sec:quasi} we 
give the definition of quasi-Poisson manifolds and of moment maps
for Hamiltonian quasi-Poisson manifolds.
Section \ref{sec:examples} describes the fundamental examples of 
$M=G$ with the action by conjugation, and of the
conjugacy classes in $G$. 
Section \ref{sec:cohomology} contains the definitions
of the {\em quasi-Poisson cohomology} and of the {\em modular class} of
a quasi-Poisson manifold. 
In Sections 
\ref{sec:fusion} and \ref{sec:reduction}, we study the {\it fusion} and 
{\it reductio
n} of quasi-Poisson manifolds. In Section \ref{sec:expo},
we show how to construct Hamiltonian quasi-Poisson manifolds 
out of Hamiltonian Poisson manifolds using the process of 
{\it exponentiation}. In Section \ref{sec:cross} we define a generalized
dynamical $r$-matrix
and we prove 
the {\it cross-section theorem}, which is our main
technical tool. Using this result, 
we show, in Section \ref{sec:foli}, that every Hamiltonian quasi-Poisson 
manifold is foliated into non-degenerate leaves.  In Section 
\ref{sec:P=qH} we
 prove the equivalence of the notion of a
non-degenerate Hamiltonian quasi-Poisson manifold with that of a
quasi-Hamiltonian manifold
in the sense of \cite{al:mom}, and we apply this result to 
describe the Poisson structure on the
representation variety, $\Hom(\pi_1(\Sig),G)/G$. 
In Appendix A
we explain the  relation between the quasi-Poisson bivector
on the group $G$ and the Poisson bivector on the dual of the
Lie algebra of the central extension of the corresponding
loop group $LG$. In Appendix B, we prove that the new $r$-matrix introduced
in Section 8 is indeed a
solution of a generalized classical dynamical Yang-Baxter equation.

\vskip.3in

\begin{center}
{\sc Notations}
\end{center}
\vskip.1in
Let $G$ be a Lie group with Lie algebra $\g$. 
Although many of the results of this paper hold for non-compact Lie
groups, we shall assume for simplicity that $G$ is compact.
For any $G$-manifold $M$ and any $\xi\in\g$, 
the generating vector field 
of the induced 
infinitesimal action is defined by 
$\xi_M(m):=\f{d}{d t}|_{_{t=0}}\exp(-t\xi). m$, 
for $m \in M$.
The Lie algebra
homomorphism $\g\to C^\infty(M;TM),\,\xi\mapsto \xi_M$, 
extends to an equivariant map, 
$$ \wedge^\bullet \g \to C^\infty(M;\wedge^\bullet TM),$$ 
preserving wedge products and Schouten brackets.\footnote{
We briefly recall the main properties of the Schouten bracket:
$$ [\alpha,\beta]=-(-1)^{(|\alpha|-1)(|\beta|-1)}[\beta,\alpha],$$
$$ [\alpha,[\beta,\gamma]]= [[\alpha,\beta],\gamma]
            +(-1)^{(|\alpha|-1)(|\beta|-1)} [\beta,[\alpha,\gamma]],$$
$$ [\alpha,\beta \wedge \gamma] 
=[\alpha,\beta] \wedge \gamma +(-1)^{(|\alpha|-1)|\beta|} 
\beta \wedge [\alpha,\gamma].$$
 For any vector field $X$, the bracket 
$[X,\alpha]$ is just the Lie derivative.} More generally, 
for any function $\alpha\in C^\infty(M,\wedge\g)$, we denote 
by $\alpha_M$ the multi-vector field, 
$ \alpha_M(m)=(\alpha(m))_M(m).$

If $P$ is a bivector on $M$, 
then $P^{\sharp}$ denotes the
induced map from differential forms to vector fields, with the
convention $P^{\sharp}(a)(b)=P(a,b),$ for $1$-forms $a$ and $b$. 

We shall denote the left- and 
right-invariant multivector fields on $G$  generated by an element
$\beta\in\wedge^\bullet\g$ by
$\beta^L$ and $\beta^R\in C^\infty(G;\wedge^\bullet TG)$ respectively.  
The vector fields $\xi^L$ for $\xi\in\g$ are the generating 
vector fields for the right action, 
$(g,m) \in G \times G \mapsto g. m=m\,g^{-1} \in G$
and $-\xi^R$ are the generating vector fields for the left action,
$(g,m) \mapsto g. m=gm$. 

Our definitions involve the choice of an invariant
inner  product (positive definite,
non-degenerate symmetric bilinear form)
on $\g$, which we 
also use to identify $\g^*$ with $\g$.
We denote the inner product by a dot.
Using the inner product, one can define the canonical invariant
skew-symmetric $3$-tensor on $\g$, sometimes called the Cartan $3$-tensor.
In terms of an orthonormal basis $(e_a)$ of $\g$, this 
canonic
al element, $\phi\in\wedge^3\g$, is given by
$$ \phi=\f{1}{12} f_{abc} e_a \wedge e_b \wedge e_c \ ,$$
where $f_{abc}=e_a \cdot [e_b,e_c]$ are 
the structure constants of $\g$. Here and below, we take the sum over
repeated indices. (Normalizations for the element $\phi$ vary in the
literature.) 

We denote by $\theta^L$ the left-invariant Maurer-Cartan 
form on $G$, and by $\theta^R$ the right-invariant Maurer-Cartan 
form. 
Let $\theta_a^L$ and $~\theta_a^R$ be the components of the 
Maurer-Cartan forms in the basis $(e_{a})$. Then 
$\iota(e_a^L)\theta_b^L=\delta_{ab}$ and 
$\iota(e_a^R)\theta_b^R=\delta_{ab}$. 
If $A:\,\g\to\g$ is a linear map, we define its 
components by $A e_b= A_{ab}e_a$. 
This gives $(\ad_\xi)_{ab}=-f_{abc}\xi_c$, for $\xi \in \g$. 
Also, at a point $g \in G$, $e_a^R=(\Ad_g)_{ab}e_b^L$,  
and $\theta_a^R=(\Ad_g)_{ab}\theta_b^L$.

\section{Quasi-Poisson manifolds}\label{sec:quasi}
Recall that a Hamiltonian Poisson $G$-manifold is a triple, 
$(M,P_0,\Phi_0)$, 
consisting of a $G$-manifold $M$, an invariant Poisson structure 
$P_0$, and an equivariant {\em moment map} $\Phi_0:\,M\to \g^*$ 
satisfying the condition, 
\begin{equation}\label{eq:Poissonmoment}
P_0^{\sharp}(\d\l\Phi_0,\xi\r)=\xi_M \ ,
\end{equation}
for all $\xi \in \g$. 
The simplest example of a Hamiltonian Poisson $G$-manifold is 
$M=\g^*$ with its linear Poisson structure and the coadjoint action; 
the identity map is then a moment map. 

\medskip

In this paper we will study a notion of Hamiltonian {\em quasi}-Poisson 
$G$-man
ifold, for which the moment map takes values in the 
group $G$ itself. The terminology ``quasi'' is motivated 
by the relation to quasi-Poisson Lie groups (see \cite{al:ma}),
which are the classical limits of quasi-Hopf algebras.
For any $G$-manifold $M$ the Cartan 3-tensor $\phi$ corresponding to 
an invariant inner product on $\g$ gives rise 
to an invariant trivector field $\phi_M$ on $M$.

\begin{definition}\label{def:quasiPoisson}
A {\em quasi-Poisson manifold} is a $G$-manifold $M$, equipped with an 
invariant bivector field $P\in C^\infty(M;\wedge^2 TM)$ such that
\begin{equation}\label{eq:Schouten}
[P,P]=\phi_M.
\end{equation}

\end{definition}

We note that if the group $G$ is Abelian, $P$ is a Poisson structure 
in the usual sense. We denote by 
$\{f,g\}=P(\d f, \d g)$ the bracket defined by $P$. It does not in general 
satisfy the Jacobi identity; instead
$$ \{\{f_1,f_2\},f_3\}+\{\{f_2,f_3\},f_1\}+\{\{f_3,f_1\},f_2\}=
2\phi_M(\d f_1,\d f_2, \d f_3).$$
An equivariant map $F:\,M\to N$ between quasi-Poisson 
$G$-manifolds $(M,P)$ and $(N,Q)$ will be called a 
{\em quasi-Poisson map} if  $\{F^* f_{1}, F^* f_{2}\}=F^*\{f_{1},f_{2}\}$ for 
all $f_{1},f_{2}\in C^\infty(N,\R)$. 
Equivalently, $F_*(P_{m})=Q_{F(m)}$ for all $m\in M$. 

To motivate the moment map condition for quasi-Poisson manifolds, 
we observe that the moment map condition for Poisson manifolds 
\eqref{eq:Poissonmoment} is equivalent to, 
\begin{equation} \label{eq:usual}
P_{0}^{\sharp} (\d(\Phi_{0}^{*}f)) = (\Phi_{0}^{*}(\ca{D}_{0}f))_{M} \ ,
\end{equation}
for all 
functions $f \in C^{\infty}(\g^{*},\R)$. 
Here $\ca{D}_0f$ is the exterior differential of $f$, viewed 
as a $\g$-valued function on $\g^*$. Equivalently, $\ca{D}_0$ is 
the $\g$-valued differential operator on $\g^*$ such that 
$(\ca{D}_0f)_a=\f{\p f}{\p \xi_a}$ in a basis 
$(e_a)$ of $\g$, defining coordinates $\xi_{a}$
on $\g^{*}$. Equation \eqref{eq:usual} becomes
\begin{equation} \label{eq:orthonormal}
P_{0}^{\sharp} (\d(\Phi_{0}^{*}f)) = \Phi_{0}^{*}
\big(\f{\p f}{\p\xi_a}\big)(e_{a})_{M}.
\end{equation}
In the non-linear case, we assume that the basis $(e_a)$
is orthonormal and we replace $\ca{D}_{0}$  by 
the $\g$-valued differential
operator $\ca{D}$ on $G$, with components
$$
({\ca{D}f})_{a}= \hh (e_{a}^{L}+e_{a}^{R})f \ .
$$ 
We note that $\ca{D}$, unlike $\ca{D}_{0}$, 
depends upon the choice of the inner
product on $\g$.

\begin{definition} \label{def:mom}
An $\Ad$-equivariant map $\Phi:\,M\to G$ is called a {\em moment map} 
for the quasi-Poisson manifold $(M,P)$ if
\begin{equation}\label{eq:pullbacks0}
P^{\sharp}(\d (\Phi^* f))=(\Phi^*(\ca{D}f))_{M} \ ,
\end{equation}
for all functions $f\in C^{\infty}(G,\R)$. 
The triple $(M,P,\Phi)$ is then called a {\em Hamiltonian quasi-Poisson
manifold}. 
\end{definition}
In the basis $(e_a)$, the moment map condition 
\eqref{eq:pullbacks0} reads 
\begin{equation}\label{eq:pullbacks}
P^{\sharp}(\d (\Phi^* f))={\f {1} {2}}\Phi^*((e_a^L+e_a^R) f)\,(e_a)_M \ .
\end{equation}
The following Lemma gives an equivalent 
version of the moment map condition. 

\begin{lemma}
Let $(M,P)$ be a quasi-Poisson $G$-manifold. An 
$\Ad$-equivariant map $\Phi:\,M\to G$ is a moment map if and only if
\begin{equation}\label{eq:moment}
P^{\sharp}(\Phi^*\theta^R_a)={\f {1}{2}} (1+\Ad_{\Phi})_{ab}(e_b)_M.
\end{equation}
\end{lemma}
\begin{proof}
First suppose that $\Phi$ satisfies \eqref{eq:moment}.
Using $\d f=(e_a^R f)\theta_a^R$, we find that
$$P^{\sharp}(\d (\Phi^* f))={\f {1} {2}} \Phi^* (e_a^R f) 
(1+\Ad_\Phi)_{ab}\,(e_b)_M \ .$$ 
Equation \eqref{eq:pullbacks} follows since $e_a^R=(\Ad_g)_{ab} e_b^L$. 
The converse is proved similarly, since, for all $g\in G$, 
one can always find $f\in C^\infty(G,\R)$ such that $\d f=
\theta_{a}^R$ at $g$. 
\end{proof}

The definition of a moment map for a quasi-Poisson manifold can also
be cast in the more invariant form,
\begin{equation}\label{eq:moment0}
P^{\sharp}(\Phi^*(\theta^R \cdot \xi))={\f {1} {2}}
((1+\Ad_{\Phi^{-1}})\xi)_M \ ,
\end{equation}
for all $\xi \in \g$,
from which it is clear that Definition \ref{def:mom} 
coincides with the definition given in
\cite{al:ma}.

\section{Examples: the Lie group and its conjugacy classes}\label{sec:examples}
The basic example of a Hamiltonian quasi-Poisson manifold is 
the group $G$ itself. Let $\psi\in\wedge^2(\g\oplus\g)$ be 
the element 
\begin{equation}\label{eq:psi} 
\psi=\hh e_a^1\wedge e_a^2,
\end{equation}
where the superscripts refer to the respective $\g$-summand. 
%Note that $\psi$ is invariant under the diagonal $\g$-action. 
A straightforward calculation
%
%\footnote{Here is the calculation in detail: \beq [e_a^1\wedge e_a^2,
%e_b^1\wedge e_b^2] &=& [e_a^1\wedge e_a^2, e_b^1]\wedge e_b^2 -
%e_b^1\wedge [e_a^1\wedge e_a^2, e_b^2]\\ &=& [e_a^1,e_b^1]\wedge
%e_a^2\wedge e_b^2 - e_b^1\wedge e_a^1\wedge [e_a^2, e_b^2]\\ &=&
%f_{abc}(e_c^1\wedge e_a^2\wedge e_b^2 +e_a^1\wedge e_b^1\wedge e_c^2)
%\eeq }
%
shows that the Schouten bracket of $\psi$ is the following element in
$\wedge^3(\g\oplus\g) \ ,$ 
$$ [\psi,\psi]=\f{1}{4}f_{abc}(e_a^1\wedge e_b^1\wedge e_c^2
+e_a^1\wedge e_b^2\wedge e_c^2).$$
In terms of the 
map, $\diag:\,\wedge^\bullet\g\to \wedge^\bullet(\g\oplus\g)$, 
induced by the diagonal embedding $\g\to \g\oplus\g$, this can 
be written
\begin{equation}\label{eq:nice}
[\psi,\psi]=\diag(\phi)-\phi^1-\phi^2 \ ,
\end{equation}
where $\phi$ is the Cartan $3$-tensor of $\g$.
Consider the map 
$\g\oplus\g\to C^\infty(G;TG), (\xi^1,\xi^2)\mapsto 
(\xi^2)^L- (\xi^1)^R$ given by the generators of the 
$G\times G$-action on $G$, where the 
first $G$-factor acts by the left action and the 
second $G$-factor by the right action.
The image of $-\psi$ under the
extended map 
$ \wedge^\bullet(\g\oplus\g)\to C^\infty(G;\wedge^\bullet TG)$
is the bivector field on $G$,
$$ P_G=\hh e_a^R\wedge e_a^L \ . $$
When the action map is composed with the diagonal
embedding of $\g$ in $\g \oplus \g$, the elements of $\g$ are mapped to
the vector fields generating the conjugation action. We denote
this infinitesimal action, and its extension to $\wedge^\bullet\g$, 
by $\xi \mapsto \xi_{G}$. Therefore, 
Equation \eqref{eq:nice}
yields
$$ [P_G,P_G]=\phi_G-\phi^L+\phi^R.$$
Since $\phi$ is $\Ad$-invariant, $\phi^L$ equals $\phi^R$ and 
therefore

$$ [P_G,P_G]= \phi_{G} \ .$$
This proves the first part of
\begin{proposition}
Let $M=G$, with $G$ acting on itself 
by conjugation, and let $P_G$ be the bivector field
\begin{equation}
\label{eq:Gbivector}
P_G=\hh e_a^R \wedge e_a^L.
\end{equation}
Then $(M,P_G)$ is a Hamiltonian quasi-Poisson $G$-manifold with 
moment map $\Phi:\,M\to G$ the identity map. 
\end{proposition}
 
\begin{proof}
It remains to verify the moment map condition \eqref{eq:pullbacks}, 
which in this case is
\begin{equation}\label{eq:toshow}
P_{G}^{\sharp}(\d f) = \hh ((e^{L}_{a}+e_{a}^{R})f) (e_{a})_G.
\end{equation}
 From the definition of $P_G$ we obtain,
$$
P_{G}^{\sharp}(\d f) = \hh ((e_{a}^{R}f) e^{L}_{a}- (e_{a}^{L}f) e_{a}^{R}).
$$
These two expressions differ by
$(e_{a}^{L}f)e_{a}^{L}- (e_{a}^{R}f)e_{a}^{R} \ ,$
which vanishes by $\Ad$-invariance of the inner product.
\end{proof}
\medskip
\begin{remark} \label{rem:LGKKS}
In Appendix A, we shall give a heuristic derivation of 
the bivector field $P_G$ 
as a quotient of a formal Poisson bivector $P_{0,L\g^*}$ on the dual
of the loop algebra $L\g$ of $\g$.
\end{remark}
\bigskip

It is well-known that for any Hamiltonian Poisson manifold, 
the moment map is a Poisson map with respect to the linear Poisson 
structure on $\g^*$. Similarly
\begin{proposition}\label{prop:moment} 
Let $(M,P)$ be a Hamiltonian quasi-Poisson $G$-manifold with moment map 
$\Phi:\,M\to G$. Then $\Phi$ is a quasi-Poisson map.
\end{proposition}

\begin{proof}
Using \eqref{eq:pullbacks} and \eqref{eq:toshow} we find that
$$
\Phi_* P^\sharp(\Phi^* \d f)= \Phi_*(\Phi^* (\ca{D} f))_M=
(\ca{D} f)_G = P_G^\sharp(\d f)
$$
for all $f \in C^\infty(G, \R)$. This shows $\Phi_* P = P_G$.
\end{proof}
\bigskip
To compare the bivector field $P
_G$ to the linear 
Poisson structure $P_{0,\g}$ on $\g\cong \g^*$, 
consider the Taylor expansion of $P_G$ near the origin, 
using the coordinates provided by the exponential map,   
$\exp:\,\g\to G$. Let $\nu(s)=\f{s}{1-e^{-s}}=1+\f{s}{2}+{\rm O}(s^2)$. 
The operator $\nu(\ad_\xi):\,\g\to\g$ is well-defined for 
$\xi$ small, and  
\begin{equation} \label{eq:nu}
e_a^L = (\nu(\ad_\xi))_{ab}~\f{\p}{\p \xi_b},\ \ \ 
   e_a^R = (\nu(\ad_\xi))_{ca}~ \f{\p}{\p \xi_c} \ ,
\end{equation}
where $\f{\p}{\p \xi_a}$ are the coordinate vector fields. Therefore, 
$$ P_G=-\hh \big((\nu(\ad_\xi))^2\big)_{ab}~\f{\p}{\p \xi_a}\wedge
\f{\p}{\p \xi_b}
=-\hh (\ad_\xi)_{ab}~\f{\p}{\p \xi_a}\wedge 
\f{\p}{\p \xi_b}+{\rm O}(||\xi||^2),$$
showing that the linearization of $P_G$ is the linear Poisson structure,
\begin{equation} \label{eq:linear}
P_{0,\g}=\hh f_{abc}\xi_c~\f{\p}{\p \xi_a}\wedge\f{\p}{\p \xi_b}=
-\hh (\ad_\xi)_{ab}~\f{\p}{\p \xi_a}\wedge\f{\p}{\p \xi_b}. 
\end{equation}

\begin{proposition} \label{PCo}
 For every conjugacy class $\Co\subset G$, the bivector field 
$P_{G}$ on $G$ is tangent to $\Co$. Thus $\Co$ is a quasi-Poisson 
manifold, with moment map the embedding $\Phi:\Co\hra G$.
\end{proposition}
\begin{proof}
By Equation \eqref{eq:toshow}, $P_G^\sharp:\,T^*G\to TG$ takes values 
in the image of the action map $\g\to TG,\ \xi\to\xi_G$. Equivalently, 
$P_G$ is tangent to the orbits for the conjugation action.  
\end{proof}
To obtain a more explicit description of the bivector field 
on a conjugacy class $\Co$, identify the tangent space 
$T_g\Co$ at $g\in \Co$ with $\g_g^\perp$, 
where $\g_g=\{\xi\in\g|\,\Ad_g\xi=\xi\}$ is the Lie algebra of the 
stabilizer of $g$. The operator $\Ad_g-1$ is invertible on 
$\g_g^\perp$. Using this inverse, we set 
$$ \Big(\f{\Ad_g+1}{\Ad_g-1}|\g_g^\perp \Big):\,\g_g^\perp\to \g_g^\perp
\ .$$
We will view this as a linear operator 
on $\g$ acting trivially on $\g_g$. 
We claim that the bivector $P_G$ at the point $g\in G$ can be written 
\begin{equation} \label{eq:equationP}
P_G =
\f{1}{4}\Big(\f{\Ad_g+1}{\Ad_g-1}|\g_g^\perp\Big)_{ab}\,(e_a)_G 
\wedge (e_b)_G \ , 
\end{equation}
showing explicitly that $P$ is tangent to the orbits. Indeed, 
\beq
P_G
&=& \f{1}{4}(\Ad_{g^{-1}}-\Ad_{g})_{ab}e_a^L\wedge  e_b^L 
= \f{1}{4}(\Ad_g+1)_{cb}(\Ad_{g^{-1}}-1)_{ac}e_a^L\wedge  e_b^L\\
&=& \f{1}{4}(\Ad_g+1)_{cb}(\Ad_g-1)_{ca}e_a^L \wedge e_b^L
= - \f{1}{4}(\Ad_g+1)_{ab}(e_a)_G\wedge \,e_b^L \ .
\eeq 
which yields \eqref{eq:equationP}, since $e^{L}_{b}= 
(((1-\Ad_{g})|\g^{\perp}_{g})^{-1})_{ab}(e_{b})_{G}$.

Notice that, in analogy to \eqref{eq:equationP}, the  
linear Poisson bivector $P_{0,\g}$ at the point $\xi \in \g$ 
can be written
\begin{equation}\label{eq:equationP0}
P_{0,\g}= \hh ((\ad_\xi|\g_\xi^\perp)^{-1})_{ab}
(e_a)_\g \wedge (e_b)_\g \ , 
\end{equation}
where $G$ acts on $\g$ by the adjoint action.

\section{Quasi-Poisson cohomology} \label{sec:cohomology}
On a Poisson manifold $(M,P)$, the graded algebra of multivectors, 
$C^\infty(M; \wedge^\bullet TM)$, 
with the differential $\d_P = [P,~ \cdot ~]$, is a complex. In fact, 
since $\d_P^2= \frac{1}{2} [[P,P], ~ \cdot ~ ]$, the property
$[P,P]=0$ ensures that $\d_P$ squares to zero. 
The cohomology
of $\d_P$ is called the Poisson cohomology of $(M,P)$.

Let $(M,P)$ be a quasi-Poisson $G$-manifold. Then, $\d_P:=[P, ~ \cdot ~ ]$
defines an operator on the space of multivectors. Its square
is in general non-vanishing, $\d_P^2=\frac{1}{2} [\phi_M, ~ \cdot ~ ]$.
However, when restricted to the subspace of $G$-invariant
multivectors, $C^\infty(M; \wedge^\bullet TM)^G$, $\d_P$ becomes
a differential.

\begin{definition}
Let $(M,P)$ be a quasi-Poisson $G$-manifold. The {\em quasi-Poisson cohomology}
of $(M,P)$ is the cohomology of the differential
$\d_P=[P, ~ \cdot ~ ]$ on the space of $G$-invariant multivectors,
$C^\infty(M; \wedge^\bullet TM)^G$.
\end{definition}

Let $M$ be an orientable manifold and let $*_\mu$ be the 
isomorphism from multivectors to
differential forms on $M$ defined by a volume form, $\mu$. 
The de Rham differential
on the space of  differential forms, $\d  \ ,$ translates into the 
operator $\partial_\mu:= - *_\mu^{-1}
\circ \d \circ *_\mu$ on multivectors, called the 
{\em Batalin-Vilkovisky} or {\em BV-operator}. 
The BV-operator, $\partial_\mu$, is a 
generator of the Schouten bracket 
on $C^\infty(M; \wedge^\bullet TM)$, that is,
$$
[\alpha, \beta] = (-1)^{|\alpha|}(\partial_{\mu} (\alpha \wedge
 \beta) -
(\partial_{\mu} \alpha)\wedge \beta - (-1)^{|\alpha|} \alpha
\wedge (\partial_{\mu} \beta)) \ ,
$$
where $|\alpha|$ is the 
multivector degree of $\alpha$.
Since, moreover, the BV-operator is of square $0$, it is  a 
(super-) derivation of the Schouten bracket,
$\partial_\mu [\alpha, \beta]=[\partial_\mu \alpha, \beta]+
(-1)^{|\alpha|-1} [\alpha, \partial_\mu \beta]$. 

\begin{lemma}
Assume that $\mu$ is a $G$-invariant volume form on a $G$-manifold $M$.
Then, the map 
$\wedge^\bullet \g \rightarrow C^\infty(M; \wedge^\bullet TM)$
induced  by the $G$-action on $M$, $\xi \mapsto \xi_M$,
is a homomorphism of complexes with respect to the Lie algebra
homology operator, $\partial: \wedge^\bullet \g \rightarrow 
\wedge^{\bullet -1} \g$ and the BV-operator, $\partial_\mu$.
\end{lemma}
\begin{proof}
The map $\xi \mapsto \xi_M$ is a homomorphism
with respect to both the exterior product and the Schouten
bracket. Moreover, the operator $\partial$ is a generator 
of the Schouten bracket on $\wedge^\bullet \g$.
Hence, it is sufficient to prove the property $\partial_\mu \xi_M=
(\partial \xi)_M$ on the elements of degree 1. 
 For all $\xi \in \g$, $\partial \xi =0$. To compute
$\partial_\mu \xi_M = 
- *^{-1}_\mu \circ \d \circ *_\mu \xi_M \ ,$ we consider
$$\d (*_\mu  \xi_M) = \d \iota(\xi_M) \mu =\L(\xi_M) \mu \ ,$$
where $\iota$ denotes an interior product and $\L$ a Lie derivation.
Since $\mu$ is $G$-invariant, $\L(\xi_M) \mu=0$, whence 
$\partial_\mu (\xi_M) =0$,
which shows that $\p _{\mu}(\xi_{M})=(\partial \xi)_M$, for all $\xi \in \g$.
\end{proof}

We  define the {\em modular vector field}
on a quasi-Poisson $G$-manifold $(M,P)$ 
with given $G$-invariant volume form $\mu$ by the formula, 
$X_\mu:= \partial_\mu P$.

\begin{proposition}
 For any $G$-invariant volume form $\mu$ on the quasi-Poisson $G$-manifold 
$(M,P)$, the modular vector field $X_\mu$ is $G$-invariant and  
a cocycle with respect to $\d_P$.
The quasi-Poisson cohomology class of $X_\mu$ 
is independent of the choice of $\mu$.
\end{proposition}

\begin{proof}
Since the BV-operator is a derivation of the Schouten bracket,
for each $\xi \in \g$, $[\xi_{M},X_{\mu}] = [\xi_{M},\p_{\mu}P]
=\p_{\mu}[\xi_{M},P]-[\p_{\mu}\xi_{M},P] = 0$. Moreover,
$\d_P X_\mu = [P, \partial_\mu P]= \frac{1}{2} \partial_\mu \phi_M$.
The element $\phi$ defines a cycle in Lie algebra homology. Hence, we obtain
$\partial_\mu \phi_M= (\partial \phi)_M=0$ and $\d_P X_\mu=0$.

Choosing $\tilde{\mu} = f \mu$ with $f$ a positive 
$G$-invariant function on $M$,
one obtains the new BV-operator, $\partial_{\tilde{\mu}} = \partial_\mu
- \iota(\d (\ln f))$.
The modular vector field also changes, 
$$
X_{\tilde{\mu}} = X_\mu - \iota(\d (\ln f)) P=
X_\mu + [P, \ln(f)] = X_\mu + \d_P(\ln f).
$$
We conclude that 
the class of $X_\mu$ in the quasi-Poisson cohomology 
is independent of the choice of $\mu$.
\end{proof}
We refer to the class of $X_\mu$ in the quasi-Poisson cohomology
as the {\em modular class} of the quasi-Poisson $G$-manifold $(M,P)$.
This definition extends the standard definition of the modular
class of a Poisson manifold \cite{we:mo}.
\begin{proposition}
The modular class of the quasi-Poisson $G$-manifold $(G,P_{G})$, where
$P_{G}$ is defined by Formula \eqref{eq:Gbivector}, vanishes.
\end{proposition}
\begin{proof}
Let $\mu$ be the bi-invariant volume form on $G$ defined by the basis
$(e_a)$ of the Lie algebra $\g$. Then, 
$X_\mu=\partial_\mu P_{G}= - *_\mu^{-1} \d *_\mu P_{G}$,
where $*_\mu P_{G}= \hh \iota(e_a^R) \iota(e_a^L) \mu$. Applying the
de Rham differential yields
\beq
\d \iota(e_a^R) \iota(e_a^L) \mu & = &
{\L}(e_a^R) \iota(e_a^L) \mu -  \iota(e_a^R) \d \iota(e_a^L) \mu \\
&=&  \iota(e_a^L) {\L}(e_a^R) \mu - \iota(e_a^R) {\L}(e_a^L) \mu =0 \ ,
\eeq
since $e_{a}^{R}$ and $e_{a}^{L}$ commute, and since
$\mu$ is both left-and right-invariant and closed. This implies 
$X_\mu=0$.
\end{proof}

\section{Fusion of quasi-Poisson manifolds}\label{sec:fusion}
Any Hamiltonian Poisson $G\times G$-manifold becomes
a Hamiltonian $G$-manifold for the diagonal $G$-action, with 
moment map the sum of the moment map components. For Hamiltonian 
quasi-Poisson manifolds a similar statement is true using the 
{\em product} of the moment map components. However, it 
is necessary to change the bivector field.
In this section $H$ is a compact Lie group with an invariant inner
product, possibly $H=\{e\}$.

\begin{proposition}\label{prop:fusion} 
Let $(M,P)$ be a quasi-Poisson $G\times G\times H$-manifold. Then 
\begin{equation}
P_{fus}:=P-\psi_M \ ,
\end{equation}
where $\psi_M$ is the image of \eqref{eq:psi} under 
the $G\times G$-action map, defines a quasi-Poisson structure 
on $M$ for the diagonal $G\times H$-action. Moreover, if 
$(\Phi^1,\Phi^2,\Psi):\,M\to G\times G\times H$ is a moment map for the 
$G\times G\times H$-action, the pointwise product 
$$(\Phi^1\Phi^2,\Psi)$$
is a moment map for the diagonal action.   
\end{proposition}
\begin{proof}
The trivector field for the $G\times G$-action is the sum  
of the trivector fields $\phi_M^1$ and $\phi_{M}^{2} \,$ for the two $G$-factors.
By assumption, $P$ is $G \times G$-invariant and $[P,P] =
\phi^{1}_{M}+\phi^{2}_{M}$. Therefore $[P,\psi_{M}]=0$, and  
$$ [P-\psi_{M},P-\psi_{M}] =
\phi^{1}_{M}+\phi^{2}_{M}+[\psi_{M},\psi_{M}] \ .$$
By Formula \eqref{eq:nice},
this implies that
$$[P,P]=(\diag(\phi))_{M} = \phi^{\diag}_{M} \ ,$$
where $\phi^{\diag}_{M}$ is the image of $\phi$ under the map
extending the infinitesimal diagonal action of $G$ on $M$.

Now suppose that the action is Hamiltonian. For any maps $\Phi^1$
and $\Phi^2$ from $M$ to $G$, 
$$ (\Phi^1\Phi^2)^*\theta_a^R
=(\Ad_{\Phi^1})_{ab}(\Phi^2)^*\theta_b^R 
+(\Phi^1)^*\theta_a^R.$$
This relation together with the moment map property,
$$P^{\sharp}(\Phi^i)^{*}\theta^{R}_{a}= {\f
{1}{2}}(1+\Ad_{\Phi^i})_{ab}(e_{b})^{i}_{M} \ , $$
for $i = 1,2$, implies,
$$
P^{\sharp}((\Phi^1\Phi^2)^*\theta_a^R) = 
 {\f {1} {2}}
(1+\Ad_{\Phi^1})_{ac}(e_c)_M^1+ {\f {1} {2}}
(\Ad_{\Phi^1})_{ab}(1+\Ad_{\Phi^2})_{bc}(e_c)_M^2 \ .
$$
By equivariance of $\Phi^{i} \ ,$
$\l (\Phi^i)^{*}\theta^{R}_{a},(e_{b})^{i}_{M} \r 
= (\Ad_{\Phi^{i}} - 1)_{ab} \ ,$ 
and therefore, 
$$\psi_{M}^{\sharp}((\Phi^1\Phi^2)^*\theta_a^R) = {\f {1}{2}}
((\Ad_{\Phi^{1}})_{ab}(1-\Ad_{\Phi^{2}})_{bc} (e_{c})^{1}_{M} 
+ {\f {1}{2}} (-1 + \Ad_{\Phi^{1}})_{ac} (e_{c})^{2}_{M} \ .   $$
We obtain 
$$P_{fus}^{\sharp}((\Phi^{1}\Phi^{2})^{*}\theta^{R}_{a})= {\f {1}{2}}
(1+\Ad_{\Phi^{1}\Phi^{2}})_{ab}((e_{b})^{1}_{M}+(e_{b})^{2}_{M}) \ ,$$
which is the moment map condition for the diagonal action.
\end{proof}
\begin{example}
Let $M=G$, viewed as a $G\times G$-space for the 
action defined by $(g_1,g_2). a = g_1 a g_2^{-1}$. 
The trivector field for this action 
is $\phi_M=\phi^L-\phi^R=0$, hence $P=0$ defines a quasi-Poisson 
structure. The diagonal action is the conjugation, and the  
quasi-Poisson structure $P_{fus}$ is defined by the bivector field, 
$P_{fus}=\hh e_a^R \wedge e_a^L \ ,$ introduced above. 
We remark that the $G\times G$-action on $M=G$ does {\em not} 
admit a moment map. 
\end{example}

\begin{example}\label{ex:double1}
Let $M=D(G):=G\times G$ with $G \times G$-action
$$(g_1,g_2). (a_1,a_2)=
(g_1 a_1 g_2^{-1}, g_2 a_2 g_1^{-1}) \ ,$$
and bivector field  
\begin{equation}  \label{eq:Pdouble}
P=\hh (e_a^{1,L}\wedge e_a^{2,R} + e_a^{1,R}\wedge e_a^{2,L}). 
\end{equation}
Then $(D(G),P)$ is a quasi-Poisson $G \times G$-manifold
obtained from 
$G\times G$, viewed as a $G^4$-space, by fusing two pairs of 
$G$-factors.  A simple calculation shows that 
$(a_1,a_2)\mapsto (a_1 a_2, a_1^{-1} a_2^{-1})$ is a moment map. 
\end{example} 

\begin{example}\label{ex:double2}
By fusing the two $G$-factors acting on the Hamiltonian 
quasi-Poisson manifold $D(G)$
of Example \ref{ex:double1}, we obtain a  
Hamiltonian quasi-Poisson $G$-manifold, 
where $G$ acts by conjugation on
each factor of $G \times G$, and the moment map is the Lie group commutator, 
$(a_1,a_2)\mapsto [a_1,a_2]=a_1 a_2 a_1^{-1} a_2^{-1}$. 
We denote this Hamiltonian quasi-Poisson $G$-manifold by ${\bf D} (G)$.
\end{example}

If $M_1$ and $M_2$ are quasi-Poisson $G$-manifolds, we denote 
their direct product with diagonal $G$-action and bivector
field $(P_1+P_2)_{fus}$ by $M_1\fus M_2$. 
We remark that the two projection mappings,
$M_1\fus M_2 \to M_j$ ($j=1,2$), are quasi-Poisson maps.

Propositions \ref{prop:fusion} and \ref{prop:moment} have the following
consequence:
\begin{corollary}
The group multiplication, $\Mult_G:\,G\fus G\to G,\ (g_1,g_2)\mapsto
g_1\,g_2 \ ,$ 
is a quasi-Poisson map.
\end{corollary}

\begin{proof}
Since the identity map is a moment map for $(G,P)$,   
Proposition \ref{prop:fusion} shows that $\Mult_G$ is a 
moment map for $G\fus G$. Hence it is a quasi-Poisson map 
by Proposition \ref{prop:moment}. 
\end{proof}
\begin{proposition}\label{prop:action}
 For any quasi-Poisson manifold $(M,P)$, the action map 
$$ \A:\, G\fus M \to M,\ \ (g,m) \mapsto g. m$$ 
is a quasi-Poisson map. 
\end{proposition}

\begin{proof}
Since $P$ is $G$-invariant, $ \A_*P=P$. 
To compute $\A_*P_G$ we observe that
$\A_*e_a^R=-(e_a)_M$, since $-e_a^R$ is the 
vector field generating the left action on $G$, and $\A$ is equivariant with
respect to this action. Thus
$$ \A_*P_G=- \hh (\Ad_g)_{ab}(e_a)_M \wedge (e_b)_M.$$ 
 Finally, since $\A_*(e_a)_M=(\Ad_{g^{-1}})_{ab}(e_b)_M$, we obtain
$$ \A_*(-\hh(e_a^L-e_a^R)\wedge  (e_a)_M)
=\hh (\Ad_{g^{-1}}-1)_{ab}(\Ad_{g^{-1}})_{ac}(e_b)_M \wedge (e_c)_M 
=\hh (\Ad_g)_{ab}(e_a)_M \wedge (e_b)_M \ ,
$$
which cancels the term $\A_*P_G$. 
\end{proof}

It is evident that the fusion operation is associative.
Given a Hamiltonian quasi-Poisson $G\times G\times G\times H$-manifold
$M$, we can begin by fusing the last two $G$-factors, or
the first two $G$-factors. The $G\times H$-equivariant bivector fields 
on $M$ thus obtained are identical. 
In the following sense, fusion is also commutative.
\begin{proposition}
Let $(M,P,(\Phi^1,\Phi^2,\Psi))$ be a Hamiltonian quasi-Poisson 
$G\times G\times H$-manifold. Let 
$$ R:\,M\to M,\ \ m\mapsto (e,\Phi^1(m),e). m$$
be the action of the second $G$-factor by the value of 
the first moment map component. 
Then $R$ is a diffeomorphism, with properties 
\begin{equation}\label{eq:Rmoment}
 R^*(\Phi^2\Phi^1)=\Phi^1\Phi^2 \
\end{equation}
and
\begin{equation}\label{eq:Rtwist}
 R_*(P-\hh 
(e_a)^{1}_M\wedge (e_a)^{2}_M)=P-\hh (e_a)^{2}_M\wedge (e_a)^{1}_M \ .
\end{equation}
\end{proposition}
\begin{proof}

Equation \eqref{eq:Rmoment} follows by equivariance of the moment map,
$$R^*(\Phi^2\Phi^1)=R^*(\Phi^2)R^*(\Phi^1)
=\Ad_{\Phi^1}(\Phi^2) ~ \Phi^1=\Phi^1  \Phi^2 \ .
$$
To prove \eqref{eq:Rtwist},
write $R$ as the composition of two maps, $R=R_2\circ R_1 \ ,$ where 
$R_1:\,M\to G\times M,m\mapsto(\Phi^{1}(m),m)$, and $R_2:\,G\times M\to M$
is the action map $\ca{A}^2$ of the second $G$-factor. The tangent map 
to $R_1$ is defined by
$$ (R_1)_*(X)=\l(\Phi^1)^*\theta_a^R,X\r e_a^R +X \ .$$
In particular, $(R_1)_*(e_a)^{1}_M=(e_a)^{1}_M 
+(\Ad_{(\Phi^1)^{-1}}-1)_{ab} e_b^R$ and 
$(R_1)_*(e_a)^{2}_M=(e_a)^{2}_M$, so that
$$ \hh (R_1)_* (e_a)^{1}_M\wedge (e_a)^{2}_M =
 \hh(e_a)^{1}_M\wedge (e_a)^{2}_M - {\f {1}{2}}(\Ad_{(\Phi^1)^{-1}}-1)_{ab} 
(e_a)^{2}_M \wedge e_b^R \ .
$$
We now use 
the fact that $(R_2)_*(e_a)^{1}_M=
(e_a)^{1}_M$ and $(R_2)_*(e_a^R)=-(e_a)^{2}_M \ ,$ and the relation
$(R_2)_*(e_a)^{2}_M =(\Ad_{g^{-1}})_{ab} (e_b)^{2}_M$. We find that
\beq
 \hh (R_2 R_1)_* (e_a)^{1}_M\wedge (e_a)^{2}_M &=&
\hh(\Ad_{(\Phi^1)^{-1}})_{ab} 
(e_a)^{1}_M\wedge (e_b)^{2}_M+{\f {1}{2}}(\Ad_{(\Phi^1)^{-1}}-1)_{ab} 
(e_a)^{2}_M \wedge (e_b)^{2}_M\\
&=&
-\hh(\Ad_{\Phi^1})_{ab} 
(e_a)^{2}_M\wedge (e_b)^{1}_M -\hh (\Ad_{\Phi^1})_{ab} 
(e_a)^{2}_M \wedge (e_b)^{2}_M.
\eeq
On the other hand, $(R_1)_*P$ is the sum of $P$, of $(\Phi^1)_*P=P_G$, 
and of cross-terms which can be found by pairing this bivector
with $\theta_a^R$ and using the moment map condition. The result is
$$ (R_1)_* P=P_G+P+{\f {1}{2}}
(1+\Ad_{\Phi^1})_{ab}~e_a^R\wedge (e_b)^{1}_M.
$$
We next apply $(R_2)_*$ to this result. 
The push-forward of 
$P$ and of $P_G$ are obtained as in the proof of 
Proposition \ref{prop:action}, and we obtain 
$$ 
(R_2 R_1)_* P
=P-\hh (\Ad_{\Phi^1})_{ab} (e_a)^{2}_M\wedge (e_b)^{2}_M
- {\f {1}{2}}(1+\Ad_{\Phi^1})_{ab}(e_a)^{2}_M\wedge
(e_b)^{1}_M.
$$ 
Equation \eqref{eq:Rtwist} is now obtained by subtracting the
expression for $\hh R_{*} (e_a)^{1}_M\wedge (e_a)^{2}_M$
from that for $R_{*}P$. 
\end{proof}

\section{Reduction of Hamiltonian quasi-Poisson manifolds}
\label{sec:reduction}
If $(M,P_0)$ is a Poisson $G$-manifold, and $M_*$ the open 
subset of $M$ on which $G$ acts freely, the quotient $M_*/G$ carries a
unique Poisson structure such that the quotient map 
is Poisson. If the action is Hamiltonian with moment 
         map $\Phi_0:\,M\to \g^*$, then $\Phi_0$ has maximal rank 
on $M_*$, and for every coadjoint orbit, $\O\subset \g^*$,
$(\Phi_0^{-1}(\O) \cap M_{*})/G$ is a smooth Poisson 
submanifold of $M_{*}/G$. 
More generally, if $M_{**}$ denotes the subset where the 
action is locally free, {\it i.e.}, the stabilizer groups are finite, 
$M_{**}/G$ is a Poisson orbifold and 
$(\Phi_0^{-1}(\O) \cap M_{**})/G$ is a Poisson suborbifold
of $M_{**}/G$. The space
$M_{\O}:= \Phi_{0}^{-1}(\O)/G$ is called the {\em reduced space}. 
 For the orbit of $0 \in \g^* \ ,$ we also use the notation 
$M\qu G:=\Phi_0^{-1}(0)/G$.

Similar assertions hold for a quasi-Poisson manifold $(M,P)$. First,
since $\phi_M$ vanishes on invariant forms, the space 
$C^\infty(M,\R)^G$ of $G$-invariant functions is a Poisson 
algebra under $\{\cdot,\cdot\}$. Therefore $M_*/G$ is a 
Poisson manifold in the usual sense. 

\begin{theorem}[quasi-Poisson reduction]
Let $(M,P,\Phi)$ be a Hamiltonian quasi-Poisson manifold. 
Then $\Phi$ has maximal rank on the subset $M_{**}$ where 
the action is locally free. For each conjugacy class 
$\Co\subset G$, the intersection of $M_\Co=\Phi^{-1}(\Co)/G$
with $M_*/G$ (resp., $M_{**}/G$) is a Poisson submanifold
(resp., suborbifold). 
\end{theorem}

\begin{proof}
Suppose that the action is locally free at the point $m \in M$.
To prove that $\Phi$ has maximal rank at $m$, we need to show 
that the equation $\Phi^*(\theta^R \cdot \xi)(m)=0$ does not 
admit a non-trivial solution $\xi\in\g$. Let $\xi\in\g$ be any 
solution of this equation. By equivariance of the moment map, 
$$
0 =\langle \Phi^* (\theta^R \cdot \xi), \eta_M \rangle (m) 
=  (\Ad_{{\Phi(m)}^{-1}} - 1) \xi \cdot \eta
$$
for all $\eta \in \g$. Hence, $\Ad_{\Phi(m)^{-1}} \xi = \xi$,  
and the moment map condition \eqref{eq:moment0} shows that
$$ 0= P_m^{\sharp}(\Phi^*(\theta^R \cdot \xi))=\xi_M(m).$$
Since the action is locally free at $m$, this implies $\xi=0$.

Because $\Phi$ is of maximal rank on $M_{**}$, 
$\Phi^{-1}(\Co)\cap M_{**}$ is a smooth submanifold 
of $M_{**}$, and $(\Phi^{-1}(\Co)\cap M_{**})/G$ is a smooth suborbifold 
of $M_{**}/G$. To show that $(\Phi^{-1}(\Co)\cap M_{**})/G$ is a Poisson 
suborbifold, we have to show that, for all invariant 
$F\in C^\infty(M,\R)^G$, the Hamiltonian vector field $P^{\sharp}(\d F)$ 
is tangent to $\Phi^{-1}(\Co)\cap M_{**}$. In fact, 
$P^{\sharp}(\d F)$ is tangent to all level sets of $\Phi$, because 
$$ (P^\sharp(\d F))\Phi^* f=-(P^\sharp(\d\Phi^* f))\,F=0$$
for all functions $f\in C^\infty(G,\R)$. Here we have used the moment
map condition \eqref{eq:pullbacks} and $(e_a)_M F=0$.
\end{proof}

 For $\Co=\{e\}$ we also write $M\qu G:=\Phi^{-1}(e)/G$.
\begin{example}
Let $\Sig$ be a compact oriented
surface of genus $h$ with $r\ge 1$ boundary
components. It is known that the representation variety 
$\Hom(\pi_1(\Sig),G)/G$ can be identified with the moduli space 
of flat $G$-connections on the trivial $G$-bundle over $\Sigma$. 
Its smooth part 
carries a natural Poisson structure
\cite{go:sy}, the leaves of which correspond to 
connections whose holonomies around the boundary circles belong to fixed
conjugacy classes in $G$. 
Using the notion of fusion, this Poisson structure can be described as
follows. First we observe that
$$  \Hom(\pi_1(\Sig),G)=\{(a_1,\ldots,a_{2h},b_1,\ldots,b_r)\in G^{2h+r}|\,
\prod_j [a_j,a_{j+1}]\prod_k b_k= e\} \ .$$
Viewing this relation 
as a moment map condition, we can write
$$ \Hom(\pi_1(\Sig),G)/G=\big(\underbrace{{\bf D}(G)\fus\ldots \fus {\bf D}(G)}_{h}\fus 
\underbrace{G\fus \ldots \fus G}_{r}\big)\qu G.$$
\end{example}

 For any Hamiltonian quasi-Poisson manifold $(M,P,\Phi)$, the embedding 
$M\to M\fus G$ defined by $m\mapsto (m,(\Phi(m))^{-1})$ is a $G$-equivariant 
diffeomorphism from $M$ onto the identity level set of the moment map for 
$M \fus G$, $(m,g) \in M\times G \mapsto \Phi(m)g \in G$. Thus it
induces a bijection from  $M/G$ to $(M\fus G)\qu G$. It is easily shown 
that on the smooth part, this bijection is a Poisson diffeomorphism. 
In particular, the Poisson structure on the representation variety 
can also be written, 
$$ \Hom(\pi_1(\Sig),G)/G=
\big(\underbrace{{\bf D}(G)\fus\ldots \fus {\bf D}(G)}_{h}\fus 
\underbrace{G\fus \ldots \fus G}_{r-1}\big)/ G.
$$
In Section \ref{sec:P=qH}, Example \ref{ex:moduli} we shall
explain that this Poisson structure coincides with
the canonical Poisson structure on the moduli space 
of flat connections on $\Sigma$.
A similar construction of the Poisson structure on the
representation variety using the properties of the 
solutions of the classical Yang-Baxter equation was developed in
\cite{fo:po}.

\section{Exponentials of Hamiltonian Poisson-manifolds}
\label{sec:expo}
Let $(M,P_0,\Phi_0)$ be a Hamiltonian Poisson $G$-manifold. We 
will show in this section that if $\Phi_0$ takes values 
in the open subset $\g^\natural \subset \g$ on which the exponential 
map has maximal rank, then $P_0$ can be modified into a 
quasi-Poisson structure with moment map $\exp \circ \, \Phi_0$. 
This construction involves a bivector field on $\g^\natural$,
$$ 
T=\hh T_{ab}e_a\wedge e_b\in C^\infty(\g^\natural ; \wedge^2 \g),
$$
where $T_{ab}(\xi)=(\varphi(\ad_\xi))_{ab}$, for $\xi \in \g^\natural$, 
and $\varphi$ is the meromorphic function of $s\in\C$,  
$$ \varphi(s)= \f{1}{s}- \hh \on{coth}
\big(\f{s}{2}\big)=-\f{s}{12}+{\rm O}(s^3).$$
We observe that $T$ is 
well-defined and smooth on $\g^\natural$, but develops 
poles on the complement of this set. 
By a result of 
Etingof and Varchenko \cite{et:ge}, $T$ is a solution of the classical
dynamical Yang-Baxter equation,
\begin{equation}\label{eq:CDYB}
 \on{Cycl}_{abc}(\f{\p T_{ab}}{\p \xi_c}+T_{ak}f_{kbl}T_{lc})=
\f{1}{4}f_{abc},
\end{equation}
where $\on{Cycl}_{abc}$ denotes the sum over cyclic permutations.  
It follows from \eqref{eq:nu} and from the relation 
$\hh (\nu(s)+\nu(-s))=1-s\varphi(s)$ 
that $T$ also satisfies the equation (see \cite{al:no})
\begin{equation}\label{eq:TAB}
T_{ab} (e_b)_\g=\f{\p}{\p \xi_a}- \hh \exp^{*}(e_a^L+e_a^R),
\end{equation}
where for a local diffeomorphism $F: M\rightarrow N$ and any
multivector field $u$ on $N$ we denote by $F^*u$ the unique 
multivector field on $M$ such that $F_*(F^*u)=u$.

%For any $\wedge\g$-valued function $\alpha \in C^\infty(M, \wedge \g)$,
%we denote by $(\alpha)_M$ the multivector field on $M$
%given by $(\alpha)_M(m) = ((\alpha)(m))_M(m)$.

\begin{theorem}\label{th:expo}
Let $(M,P_0)$ be a Hamiltonian 
Poisson $G$-manifold with moment map $\Phi_0$. 
Suppose that $\Phi_0(M)\subset \g^\natural$. Then, 
$\exp:\,\g\to G$ has maximal rank on the image of $\Phi_0$ and
$ P=P_0-(\Phi_0^*T)_M$ 
defines a Hamiltonian quasi-Poisson structure on $M$, with moment map
$\Phi=\exp \circ \, \Phi_0$.
\end{theorem}
\begin{proof}
Let $T_M=(\Phi_0^*T)_M.$ We want to compute 
$$ [P,P]=[P_0,P_0]-2 [P_0,T_M]+[T_M,T_M]. $$
The first term vanishes since $P_0$ is a Poisson structure. 
 For any function $f\in C^\infty(\g,\R)$, 
$[P_0,\Phi_0^*f]=-P_0^\sharp(\Phi_0^*\d f)$. 
Using the moment map condition \eqref{eq:orthonormal} of 
$\Phi_0$ and the invariance property $[P_0,(e_a)_M]=0$ of $P_0$,
we obtain  
\begin{equation}\label{eq:1}
[P_0,T_M ] =- \hh \Phi_0^* \f{\p T_{ab}}{\p \xi_c}
(e_a)_M \wedge (e_b)_M\wedge   (e_c)_M \, . 
\end{equation}
To calculate $[T_M,T_M]$ we use the relation 
$$ [(e_a)_M,\Phi_0^* T_{kl}]=\Phi_0^*(f_{akm}T_{ml}+f_{alm}T_{km}).$$
Taking symmetries into account, we obtain 
\beq
\lefteqn{[T_M,T_M]=
\f{1}{4}\Phi_0^*(T_{ab}T_{kl})\ [(e_a)_M\wedge  (e_b)_M , 
(e_k)_M \wedge (e_l)_M]} \\
&& \hspace{5cm}
+\Phi_0^*T_{ab} [(e_b)_M ,\Phi_0^* T_{kl}] (e_a)_M \wedge  
(e_k)_M\wedge  (e_l)_M
\\ 
&=& \Phi_0^*(T_{ab} T_{kl}) f_{bkm} (e_a)_M\wedge  (e_m)_M \wedge (e_l)_M 
+ 2 \Phi_0^* (T_{ab} f_{bkm} T_{ml}) (e_a)_M \wedge (e_k)_M\wedge  (e_l)_M 
\\
&=& \Phi_0^*(T_{ak}f_{kbl}T_{lc}) (e_a)_M\wedge  (e_b)_M 
\wedge  (e_c)_M \ .
\eeq 
Together with \eqref{eq:1}, and using the classical dynamical
Yang-Baxter equation \eqref{eq:CDYB},
this yields
$$ [P,P]=\Phi_0^*\Big(\f{\p T_{ab}}{\p\xi_c}+T_{ak}f_{kbl}T_{lc}\Big) 
(e_a)_M\wedge  (e_b)_M\wedge   (e_c)_M =\phi_M.
$$
To prove that $\Phi=\exp \circ \, \Phi_0$ is a moment map, we use 
Equation
\eqref{eq:TAB}. For all functions $f\in C^\infty(G)$,  
\beq
T_M^\sharp(\d\,\Phi^*f)&=&
\Phi_0^*T_{ab}\,(e_a)_M(\Phi_0^*\exp^*f)\ (e_b)_M  \\
&=&\Phi_0^*(T_{ab}(e_a)_\g(\exp^*f))\ (e_b)_M  \\
&=& \Phi_0^*\big(\f{\p }{\p \xi_a} (\exp^*f) 
 - \hh \exp^{*}((e_a^L+e_a^R)f)\big) \ (e_a)_M 
\\
&=& \Phi_0^*(\f{\p }{\p \xi_a}(\exp^* f))
\ (e_a)_M
- \hh \Phi^*((e_a^L+e_a^R)f) (e_{a})_{M}.
\eeq
The moment map property \eqref{eq:pullbacks} 
of $\Phi$ follows from this relation together with 
the moment map map property \eqref{eq:orthonormal} of 
$\Phi_0$.
\end{proof}
\begin{example} 
Let $P_{G}$ be the quasi-Poisson structure on $G$ defined by
\eqref{eq:Gbivector}. Under the exponential map, it pulls back to a
bivector field, $\exp^*P_{G}$, on $\g^\natural$. On the other hand, 
let $P_{0,\g}$
be the linear Poisson structure on $\g^*\cong\g$. Then 
$\exp^*P_{G}=P_{0,\g}-T_\g$.
This follows from Equation \eqref{eq:TAB}, together with the
expressions \eqref{eq:equationP} for $P_{G}$ and \eqref{eq:equationP0} for
$P_{0,\g}$. 
\end{example}
\medskip
The process described in Theorem \ref{th:expo} can also be reversed.  
\begin{corollary}\label{cor:log}
Let $(M,P,\Phi)$ be a Hamiltonian quasi-Poisson $G$-manifold. Suppose
that the image of $\Phi$ is contained in an open subset $U\subset G$ 
on which the exponential map admits a smooth inverse, 
$\log:\,U\to \g$. Set $\Phi_0=\log \circ \, \Phi$, and 
$P_0=P+T_M$. Then $(M,P_0,\Phi_0)$ is a Hamiltonian Poisson $G$-manifold
in the usual sense.   
\end{corollary}

\section{Cross-section theorem}
\label{sec:cross}
The Guillemin-Sternberg symplectic cross-section theorem states 
that for any symplectic $G$-manifold with moment map $\Phi_0:\,M\to
\g^*$, the inverse image $Y=\Phi_0^{-1}(U)$ of a slice $U\subset \g^*$ 
at a given point $\xi \in\g^*$ is a $G_\xi$-invariant 
symplectic submanifold, with moment map the restriction of 
$\Phi_0$. Thus, $Y$ is a Hamiltonian $G_\xi$-manifold, called a
{\it symplectic cross-section} of $M$.
In this section we will prove a similar result for the 
quasi-Poisson case. 

Given $g\in G$, we denote the stabilizer
of $g$ with respect to the conjugation action of $G$
on itself by $H =G_g$. Any sufficiently small
connected  open neighborhood $U$ of $g$ in $H$ is a 
slice for this action. That is, the natural map 
$ G\times_H U \to G.U$, 
$(g,h) \mapsto ghg^{-1}$, is a diffeomorphism onto its image.
There is an orthogonal decomposition 
\begin{equation}\label{eq:splitu}
TG|_U=TU\oplus (U\times\h^\perp),
\end{equation}
where $\h={\rm Lie}(H)$ and
the second summand is embedded by means of the 
vector fields generating the adjoint action.

Let $(M,P,\Phi)$ be a Hamiltonian quasi-Poisson manifold. 
By equivariance of $\Phi$, the inverse image $Y=\Phi^{-1}(U)$ is a smooth 
$H$-invariant submanifold of $M$, and there is an $H$-equivariant 
splitting of the tangent bundle,
\begin{equation}\label{eq:splitting}
TM|_Y=TY\oplus (Y\times \h^\perp).\end{equation}
Dually, 
\begin{equation}\label{eq:splitting1}
T^*M|_Y=T^*Y\oplus (Y\times (\h^\perp)^*) \ .
\end{equation}
\begin{lemma}\label{lem:spltting}
The splitting \eqref{eq:splitting1} is $P|_Y$-orthogonal, that is,   
there is a decomposition, 
$P|_Y=P_Y+P_Y^\perp$, where $P_Y$ is a 
bivector field on $Y$ and $P_Y^\perp\in 
C^\infty(Y ; \wedge^2\h^\perp)$. 
\end{lemma}
\begin{proof}
At $m\in Y \subset M$,
the fiber $\{m\}\times (\h^\perp)^*$ is the space of
covectors $\alpha=\Phi^*(\theta^R \cdot \xi)(m)$ with 
$\xi\in\h^\perp$. Let $\xi'=\hh (1+\Ad_{\Phi(m)^{-1}})\xi\in\h^\perp$. 
By the moment map condition \eqref{eq:moment0},
it follows that, 
$$
 P_m(\alpha,\beta) = \l \xi'_M(m), \beta \r =0 \ ,
$$
for all $\beta\in T^*_mY$.
\end{proof}
An explicit description of the bivector $P_Y^\perp$ can be given 
in terms of the moment map. Define $r\in C^\infty(U;\wedge^2 \h^\perp )$ 
as 
\begin{equation} \label{eq:r}
r(h)=\hh r_{ab}(h) e_a \wedge e_b \ ,
\end{equation}
where 
\begin{equation} \label{eq:rmatrix} 
r_{ab}(h)= -\f{1}{2}
\Big(\f{\Ad_h+1}{\Ad_h-1}\Big|\h^\perp\Big)_{ab} \ . 
\end{equation}
Since $\Phi$ is a quasi-Poisson map,
the description \eqref{eq:equationP} 
of the quasi-Poisson structure on $G$ shows that
$P_Y^\perp=- \Phi_Y^* r$, where $\Phi_{Y}=\Phi|_{Y}$.

We will need the following property of the tensor field $r$.
 For every vector field $X\in C^\infty(U;TG)$ on $U$, let 
$X_\h$ be the orthogonal projection of $X$ to a vector field tangent to 
$U$, using 
the splitting \eqref{eq:splitu}. Let $f_{abc}^\h$ be the structure 
constants of $\h$, viewed as the components of 
a skew-symmetric tensor of $\g$ 
under the inclusion $\wedge^3\h\to \wedge^3\g$ .  
\begin{lemma}\label{lem:YB}
The tensor field $r$ defined by Equation \eqref{eq:rmatrix}
satisfies
\begin{equation}\label{eq:EV2}
\on{Cycl}_{abc}\Big(\hh(e_a^L+e_a^R)_\h r_{bc}
+ r_{ak}f_{kbl}r_{lc} \Big)=\f{1}{4}(f_{abc}-f^\h_{abc}).
\end{equation}
\end{lemma}
This lemma, which is proven in Appendix B, generalizes 
the result proved by
Etingof and Varchenko in 
\cite{et:ge} (see also \cite{al:no}) in 
the special case where $H$ is Abelian, where this equation reduces to the
classical dynamical Yang-Baxter equation.
\begin{theorem}
 For any Hamiltonian quasi-Poisson $G$-space $(M,P,\Phi)$, and 
any slice $U$ at a given $g\in G$ with stabilizer $G_g=H$, the 
cross-section $Y=\Phi^{-1}(U)$ 
is a Hamiltonian quasi-Poisson $H$-manifold, 
with bivector field  $P_Y=P|_Y-P_Y^\perp$ and moment map
$\Phi_Y=\Phi|_Y$.
Conversely, given a  Hamiltonian quasi-Poisson $H$-manifold 
$(Y,P_Y,\Phi_Y)$, let the associated bundle 
$G\times_H Y$ be equipped with the unique equivariant 
map $\Phi:\,G\times_H Y \to G$ and the unique invariant 
bivector field $P$, such that $\Phi$ restricts to $\Phi_Y$ 
and $P$ to $P_Y-\Phi_Y^*r$. Then $(G\times_H Y,P,\Phi)$ is 
a Hamiltonian quasi-Poisson $G$-manifold. 
\end{theorem}

\begin{proof}
Replacing $M$ by $G. Y$, we may assume that $M=G\times_H Y$. 
It is clear that the moment map condition for $(Y,P_Y,\Phi_Y)$ 
is equivalent to that for $(M,P,\Phi)$. We have to show that 
$ [P,P]=\phi_M$ if and only if $[P_Y,P_Y]=\phi^\h_Y \ .$
Let $\wt{P}_Y$ denote the extension of $P_Y$ to a $G$-invariant 
bivector field on $M$. Also, let 
$\wt{r}\in C^\infty(G. U,\wedge^2\g)$ be the $G$-invariant extension 
of $r$. By definition, $P=\wt{P}_Y- \Phi^*\wt{r}$, so what we need to 
show is that
\begin{equation}\label{eq:needtoshow}
([\Phi^*\wt{r},\Phi^*\wt{r}]- 2[\wt{P}_Y, \Phi^*\wt{r}])|_Y
=\phi_M|_Y-\phi^\h_Y.
\end{equation}
By a calculation similar to that of the term $[T_M,T_M]$ in the 
proof of Theorem \ref{th:expo}, 
$$ 
[\Phi^*\wt{r},\Phi^*\wt{r}]=
\Phi^*(\wt{r}_{ak} f_{kbl} \wt{r}_{lc}) 
(e_a)_M\wedge (e_b)_M \wedge (e_c)_M.
$$
To compute the term $[\wt{P}_Y,\Phi^*\wt{r}]$, first observe that 
$[\wt{P}_Y,(e_a)_M]=0$ by $G$-invariance of $\wt{P}_Y$. 
Using \eqref{eq:pullbacks}, we obtain 
$$ [P_Y,\Phi_Y^* r_{bc}]=- P_Y^{\sharp}(\d \Phi_Y^* r_{bc})
=-\Phi^*\Big(\hh(e_a^L+e_a^R)_\h r_{bc}\Big) (e_a)_Y. 
$$
Therefore, 
$$ 
([\Phi^*\wt{r},\Phi^*\wt{r}]- 2[\wt{P}_Y, \Phi^*\wt{r}])|_Y
=\Phi^*\big(\hh(e_a^L+e_a^R)_\h r_{bc}
+ r_{ak} f_{kbl} r_{lc} \Big) 
(e_a)_M \wedge (e_b)_M \wedge (e_c)_M \, .
$$
Together with Lemma \ref{lem:YB}, this proves \eqref{eq:needtoshow}.
\end{proof}

Using the construction from Section \ref{sec:expo}, 
the cross-section $Y$ can be equipped with an ordinary Poisson 
structure as follows. Since $g$ is in the center of $H=G_g$, 
$\Phi_Y'=g^{-1}\Phi_Y$ is also a moment map for $(Y,P_Y)$. For
$U$ small enough, 
the exponential map for $\h$ admits 
an inverse on $g^{-1}U\subset H$, $\log:\,g^{-1}U\to\h$.
Let 
$$\Phi_{0,Y}=\log(g^{-1}\Phi_Y),$$
and set $P_{0,Y}=P_Y+T_Y$. By Corollary \ref{cor:log}, 
$(Y,P_{0,Y},\Phi_{0,Y})$ is a Hamiltonian Poisson 
$H$-manifold.

\section{The generalized foliation of a Hamiltonian 
quasi-Poisson manifold}
\label{sec:foli}
It is a well-known result of Lie, for the constant rank case, and of 
Kirillov, for the general case, that for any Poisson-manifold 
$(M,P_0)$, the generalized distribution\footnote{A (differentiable) 
generalized distribution on a manifold $M$ is a family 
of subspaces $\D_m\subset T_mM$, such that 
for all $m\in M$, there exists a finite number of  
vector fields $X_1,\ldots,X_k\in C^\infty(M;TM)$ taking values 
in $\D = \cup_{m\in M}\D_m$ 
and spanning $\D_m$ at $m$. An in-depth discussion of generalized 
distributions can be found in Vaisman's book \cite{va:po}.} 
$\D_0=\im(P_0^\sharp)$ is 
integrable.
% where $P_0^\sharp:\,T^*M\to TM$ is the vector bundle 
% map defined by $P_0$. 

In this section we show that 
every Hamiltonian quasi-Poisson manifold 
is foliated by 
non-degenerate quasi-Poisson submanifolds.
Given a Hamiltonian quasi-Poisson manifold $(M,P,\Phi)$, 
define a generalized distribution $\D$ on $M$ by 
$$ \D_m:=\on{im}(P^\sharp_m)+T_m(G. m). $$ 
Because $G$ is compact, $\g = \im(1+\Ad_{g}) \oplus \ker(1+\Ad_{g})$,
for any $ g \in G$. By the moment map condition, 
the image of $P^\sharp_m$ 
always contains all $\eta_M(m)$ with $\eta$ in the image of the 
operator $1+\Ad_{\Phi(m)}$ on $\g$, therefore $\D_{m}$ can be
re-written
$$ \D_m = \on{\im}(P^\sharp_m)\oplus\{\xi_M(m)|(1+\Ad_{\Phi(m)})\xi=0\} $$ 
In particular, if $m \in M$ is such that $1+\Ad_{\Phi(m)}$ 
is invertible, then $\D_m=\on{im}(P^\sharp_m)$. 
\begin{definition}
\label{def:non-deg}
A Hamiltonian quasi-Poisson manifold $(M,P,\Phi)$ is non-degenerate 
if  $\D_m=T_mM$ for all $m\in M$. 
\end{definition}
 For example, the conjugacy classes
$\Co$ of a group $G$ are non-degenerate Hamiltonian 
quasi-Poisson manifolds. The decomposition of $G$ into conjugacy
classes is a special case of the following:
\begin{theorem}\label{th:foliation}
The distribution $\D$ is integrable, that is, through every 
point $m\in M$ there passes a unique maximal connected submanifold 
$N$ of $M$ such that $TN=\D|_N$. Each submanifold $N$ is 
$G$-invariant, and the restrictions of $P$ and $\Phi$ to $N$
define a non-degenerate Hamiltonian quasi-Poisson structure on $N$.
\end{theorem}

\begin{proof}
We show that $\D$ is integrable near any given point $m\in M$.
Let $g=\Phi(m)$, and $U\subset H=G_g$ be an $H$-invariant 
slice through $g$, and $Y=\Phi^{-1}(U)$ the corresponding
cross-section. Since $G. Y=G\times_H Y$, it suffices to 
show that the distribution $\D_Y$ induced by the quasi-Poisson
structure $P_Y$ is integrable. However, as explained in the previous 
section, $P_Y=P_{0,Y}-T_Y$, where $P_{0,Y}$ 
is an $H$-invariant
Poisson structure in the usual sense, and $T_Y$ is a bivector
field taking v
alues in the second exterior power of the 
$H$-orbit directions. Moreover, $P_{0,Y}$ admits a moment map
$\Phi_{0,Y}=\log(g^{-1}\Phi_Y)$, which implies that 
the image of $P_{0,Y}^\sharp$ contains
the $H$-orbit directions. Hence $\D_Y$ is {\em equal}
to the distribution defined by the Poisson-structure $P_{0,Y}$, 
and therefore integrable by the theorem of Lie and Kirillov.  
\end{proof}

\section{Non-degenerate quasi-Poisson manifolds}
\label{sec:P=qH}
Any non-degenerate bivector field $P_0$ on a manifold $M$ 
determines a non-degenerate 2-form $\om_0$ by the condition 
$\om_0^\flat=(P_0^\sharp)^{-1}$. It is well-known that under 
this correspondence, the Poisson condition $[P_0,P_0]=0$ is 
equivalent to the closure of the $2$-form $\om_0$.
 
In this section, we extend this correspondence between non-degenerate
Poisson manifolds and symplectic manifolds to the ``quasi''  case. 
While the role of the non-degenerate
Poisson manifolds is played by the non-degenerate 
%(minimally degenerate)
Hamiltonian quasi-Poisson $G$-manifolds, 
that of the symplectic manifolds is played by
the ``quasi-Hamiltonian
$G$-spaces'', introduced in \cite{al:mom}.
 First, we recall their definition, which includes the non-degeneracy
%(minimal degeneracy) 
assumption (c) below.
Then, we show that every non-degenerate
Hamiltonian quasi-Poisson manifold $(M,P,\Phi)$ carries a canonically
determined 2-form $\om$ such that $(M,\om, \Phi)$ is
a  quasi-Hamiltonian
$G$-space.

Let $\eta\in\Om^3(G)$ be
the bi-invariant closed 3-form,
$$ \eta
=\f{1}{12} f_{abc} \theta_a^R\wedge \theta_b^R\wedge \theta_c^R \ .$$
\begin{definition}\cite{al:mom}
\label{def:qham}
A quasi-Hamiltonian $G$-space is a triple $(M,\om,\Phi)$ 
where $M$ is a $G$-manifold, $\om$ an invariant differential 2-form, 
and $\Phi : M \to G$ is an $\Ad$-equivariant map, such that
\begin{enumerate}
\item 
$\d\om=\Phi^*\eta \ ,$
\item
$\iota((e_a)_M)\om=\hh
\Phi^* (\theta_a^L+\theta_a^R) \ ,$
\item
for all $m\in M$, the kernel of $\om_m$ is the space of all 
$\xi_M(m)$ such that $\xi$ is in the kernel of 
$1+\Ad_{\Phi(m)} \ .$ 
\end{enumerate}
\end{definition}
We will need the following two results that were proved in 
\cite{al:mom}. First, 
there is an exponentiation construction. Given a 
Hamiltonian symplectic $G$-manifold $(M,\om_0,\Phi_0)$ such that 
$\Phi_0(M)$ is contained in the set $\g^\natural\subset \g$
of regular values of the exponential map, $\exp$, 
one obtains a quasi-Hamiltonian $G$-space $(M,\om,\Phi)$ by setting 
$\Phi=\exp \circ \, \Phi_0$ and 
$$ \om=\om_0+\Phi_0^*\varpi,$$
where $\varpi\in\Om^2(\g)$ is the image of the closed 3-form 
$\exp^*\eta$ under the homotopy operator $\Om^\star(\g)\to
\Om^{\star-1}(\g)$. 
Secondly, there is a cross-section theorem.
Suppose that 
$(M,\om,\Phi)$ is a quasi-Hamiltonian $G$-space, and that $U\subset G$ 
is a slice at $g\in G$. Then the cross-section $Y=\Phi^{-1}(U)$ with 
the 2-form $\om_Y$ and the moment map $\Phi_Y$, defined as the 
pull-backs of $\om$ and $\Phi$, is a quasi-Hamiltonian 
$H$-space, where $H=G_{g}$. The canonical 
decomposition $TM|_Y=TY\oplus (Y\times \h^\perp)$ is $\om$-orthogonal.
Conversely, given a 
quasi-Hamiltonian $H$-space $(Y,\om_Y,\Phi_Y)$, the associated bundle
$M=G\times_H Y$ carries a unique structure of quasi-Hamiltonian 
$G$-space $(M,\om,\Phi)$ such that $\om$ and $\Phi$ pull back to
$\om_Y$ and $\Phi_Y$. 

If $(M,P,\Phi)$ is a non-degenerate Hamiltonian quasi-Poisson manifold 
such that $\Phi$ admits a smooth logarithm $\Phi_0:\,M\to \g$, 
one can define a 2-form $\om$ on $M$ in the following
way. The bivector 
$P_0=P+(\Phi_0^*T)_M$ is invertible, we denote its inverse by $\om_0 \ ,$
and we set 
$\om=\om_0+\Phi_0^*\varpi$. The following Lemma describes the relation 
between $\om$ and $P$.

\begin{lemma}\label{lem:moti}
Let $(M,P_0,\Phi_0)$ be a non-degenerate Hamiltonian Poisson
manifold, and let $\om_0$ be the symplectic form corresponding to $P_0$.
Suppose that $\Phi_0(M)$ is contained in the set of regular values of 
$\exp:\g\to G$, and let 
$\Phi=\exp \circ \, \Phi_0$, $\om=\om_0+\Phi_0^*\varpi$, and 
$P=P_0-(\Phi_0^*T)_M$. Then 
\begin{equation}\label{eq:defining}
P^\sharp\circ \om^\flat=\on{Id}_{TM}-\f{1}{4}(e_a)_M\otimes
\Phi^*(\theta_a^L-\theta_a^R).
\end{equation}
\end{lemma}
\begin{proof}
Given $m\in M$, let $U$ be a slice through $\Phi_0(m)$, 
and let  $Y=\Phi_0^{-1}(U)$ be 
the corresponding cross-section. We first evaluate both sides of 
\eqref{eq:defining} on elements of $T_mY$, and then on orbit directions. 
The 2-form $\Phi_0^*\varpi$ vanishes on $T_mY$, and the bivector 
$(\Phi_0^*T)_M$ vanishes on $T_m^*Y$. Hence 
$P^\sharp\circ \om^\flat|_{T_mY}=
P_0^\sharp\circ \om_0^\flat|_{T_mY}=\Id_{T_mY}$, which 
agrees with the right-hand side of \eqref{eq:defining}
since $\Phi^*(\theta_a^L-\theta_a^R)$ also vanishes on $T_mY$. 
To complete the proof we evaluate both sides of \eqref{eq:defining} on 
orbit directions. The moment map properties of $P$ and $\om$ yield
\beq 
P^\sharp (\om^\flat (e_a)_M)&=&
\hh (\Ad_{\Phi^{-1}}+1)_{ab} P^\sharp(\Phi^*\theta_b^R)\\
&=&\f{1}{4} (\Ad_{\Phi^{-1}}+1)_{ab} (\Ad_\Phi+1)_{bc} (e_c)_M\\
&=&\f{1}{4} (2+\Ad_{\Phi^{-1}}+\Ad_\Phi)_{ab} (e_b)_M \, , 
\eeq
and the same result is obtained by applying the right-hand side of 
\eqref{eq:defining} to $(e_a)_M$.
\end{proof}

Generalizing the Lemma, we can state the main result of this 
section:
\begin{theorem}\label{th:equivalence}
Every non-degenerate Hamiltonian quasi-Poisson manifold
$(M,P,\Phi)$ carries a unique 2-form $\om$ such that 
$(M,\om,\Phi)$ is a quasi-Hamiltonian $G$-space, and 
such that $\omega$ and $P$ satisfy Equation \eqref{eq:defining}.
Conversely, on every quasi-Hamiltonian $G$-space $(M,\om,\Phi)$ 
there is a unique bivector field $P$ such that $(M,P,\Phi)$ 
is a non-degenerate Hamiltonian quasi-Poisson $G$-manifold, 
and Equation \eqref{eq:defining} is satisfied.
\end{theorem}

\begin{proof}
Given $(M,P,\Phi)$, let $Y$ be a cross-section at $m$, as in the proof
of Lemma \ref{lem:moti}, and let $g=\Phi(m)$. 
Thus $(Y,P_Y,\Phi_Y)$ is a non-degenerate 
Hamiltonian quasi-Poisson $H$-manifold that corresponds to 
a quasi-Hamiltonian $H$-space, $(Y,\om_Y,\Phi_Y)$. Let 
$\om$ be the unique 2-form on $G. Y\subset M$ such that 
$(G. Y,\om,\Phi)$ is a quasi-Hamiltonian $G$-space, and 
$\om_Y$ is the pull-back of $\om$. We have to show that 
$\om$ satisfies Equation \eqref{eq:defining}. For orbit directions, 
this follows from the moment map conditions (see the proof of 
Lemma \ref{lem:moti}), while, for directions tangent to $Y$, 
it follows by applying the Lemma to $(Y,P_Y,g^{-1}\Phi_Y)$. 
Uniqueness is clear since the equation, 
\begin{equation} \label{eq:defining'}
\om^\flat\circ P^\sharp=\on{Id}_{T^*M}-\f{1}{4}
\Phi^*(\theta_a^L-\theta_a^R)\otimes (e_a)_M \ ,
\end{equation}
the transpose of \eqref{eq:defining}, defines $\om^\flat$ 
on the image of $P^\sharp$, while the moment map condition determines 
$\om^\flat$ on orbit directions.

The converse is proved similarly, using the cross-section theorem 
for Hamiltonian quasi-Poisson manifolds.
\end{proof}

By Theorem \ref{th:equivalence}, all the constructions and examples
for quasi-Hamiltonian $G$-spaces given in \cite{al:mom} can be translated 
into the quasi-Poisson picture. 
\begin{example} \label{ex:Pconj}
Let $\Co \subset G$ be the conjugacy class of a point $g$.
By Proposition \ref{PCo},
there exists a unique bivector $P$ on $\Co$ such that $(\Co, P, \Phi)$,
with $\Phi: \Co \rightarrow G$ the embedding, is a Hamiltonian
quasi-Poisson space. Similarly, by Proposition 3.1 of \cite{al:mom},
there exists a unique 2-form $\omega$ on $\Co$ such that
$(\Co, \omega, \Phi)$ is a quasi-Hamiltonian space.
Theorem \ref{th:equivalence} implies 
that the bivector $P$ and the 2-form $\omega$ are related by
Equation \eqref{eq:defining}.
\end{example}
\begin{example} \label{ex:Pdouble}
Let $D(G)=G\times G$ with bivector $P$ given by formula \eqref{eq:Pdouble},
and with the 2-form $\omega$ defined in Section 3.2 of \cite{al:mom},
$$
\omega=-  \frac{1}{2} 
(\theta^{1,L}_a \wedge \theta^{2,R}_a + \theta^{1,R}_a \wedge \theta^{2,L}_a).
$$
Both the left- and the right-hand 
sides of \eqref{eq:defining'} yield the same expression,
\beq
\Id_{TD(G)} & - & \frac{1}{4}\left(  
(\Ad_{a_2})_{ab} \, e^{1,L}_a \otimes \theta^{1,R}_b - 
(\Ad_{a_2^{-1}})_{ab} \, e^{1,R}_a \otimes \theta^{1,L}_b \right. \\
 & - & \left.
(\Ad_{a_1})_{ab} \, e^{2,L}_a \otimes \theta^{2,R}_b -
(\Ad_{a_1^{-1}})_{ab} \, e^{2,R}_a \otimes \theta^{2,L}_b \right) \, .
\eeq
Thus, the quasi-Poisson and quasi-Hamiltonian definitions of $D(G)$ agree.
\end{example}

In both the quasi-Poisson and the quasi-Hamiltonian settings there is a notion of reduction.
We show that these notions agree as well.

\begin{proposition} \label{prop:redeq}
Let $(M,P,\Phi)$ be a non-degenerate quasi-Poisson manifold and let
$\omega$
be the corresponding 2-form on $M$. Then, for any conjugacy class
$\Co \subset G$, the intersection of the reduced
space $M_\Co$ with $M_{*}/G$ carries a Poisson bivector 
$P^\Co$
induced by $P$  and a symplectic form $\omega^\Co$ induced by $\omega$,
such that $(P^\Co)^\sharp \circ (\omega^\Co)^\flat = \Id$.
\end{proposition}

\begin{proof}
Choose $g \in \Co$ and a slice $U$ containing $g$ and let $Y = \Phi^{-1}(U)$
be the cross-section. We observe that in both the
quasi-Poisson and the quasi-Hamiltonian
settings, the reduction of $M$ at $\Co$ is canonically isomorphic to the
reduction of $Y$ at the group unit
of $H = G_g$. The cross-section $Y$ carries the Poisson bivector
$P_{0,Y}$ and the symplectic form $\omega_{0,Y}$. 
According to Lemma \ref{lem:moti},
$(P_{0,Y})^\sharp \circ (\omega_{0,Y})^\flat = \Id_{TY}$. Hence, the same 
relation holds
for the reduced space.
\end{proof}

Next, we will show that 
the fusion operation for quasi-Hamiltonian $G$-spaces given in \cite{al:mom}
coincides with the fusion operation for quasi-Poisson spaces. 
Let 
$(M,\om,(\Phi_1,\Phi_2))$ be a quasi-Hamiltonian $G\times G$-space, 
and let $(M,P,(\Phi_1,\Phi_2))$ be the corresponding Hamiltonian 
quasi-Poisson $G\times G$-space. By \cite[Theorem 6.1]{al:mom}, 
the space $M$ with the diagonal $G$-action, the pointwise product of 
the moment map components $\Phi=\Phi_1\Phi_2$, and the fusion 2-form 
$$ \om_{fus}=\om-\hh \Phi_1^*\theta^L_a \wedge \Phi_2^*\theta^R_a \ ,$$
is a quasi-Hamiltonian $G$-space. On the other hand, Proposition 
\ref{prop:fusion} yields a bivector field $P_{fus}=P-\psi_M\ ,$ which,
together with the diagonal $G$-action and the moment map $\Phi_1\Phi_2$,
defines the structure of a Hamiltonian quasi-Poisson $G$-manifold
on $M$. 
The following proposition asserts that, as expected, 
$\om_{fus}$ corresponds to $P_{fus}$ under the equivalence 
established in Theorem \ref{th:equivalence}. 
\begin{proposition} \label{prop:fusioneq}
The bundle maps $P_{fus}^\sharp:\,T^*M\to TM$ and 
$\omega_{fus}^\flat:\,TM\to T^*M$ are related by 
\begin{equation} \label{eq:Pomfusion}
P_{fus}^\sharp \circ \omega_{fus}^\flat =
\Id_{TM} - \frac{1}{4} (e_a)_M \otimes \Phi^*(\theta^L_a - \theta^R_a),
\end{equation}
where $(e_a)_M = (e_a)^{1}_M + (e_a)^{2}_M$ are the 
vector fields generating the
diagonal $G$-action on $M$. 
\end{proposition}

\begin{proof}
We have to compute
$$
P_{fus}^\sharp \circ \omega_{fus}^\flat =
\Big(P - \frac{1}{2} (e_a)^{1}_M \wedge (e_a)^{2}_M\Big)^\sharp 
\circ 
\Big(\omega - \frac{1}{2} \Phi^*_1 \theta^L_a \wedge \Phi^*_2 \theta^R_a\Big)^\flat.
$$
We compute the four terms in the expansion of the right-hand side. 
By Theorem \ref{th:equivalence}, the first term is 
$$
P^\sharp \circ \omega^\flat =
\on{Id}_{TM}-\f{1}{4}
(e_a)^{1}_M \otimes
\Phi_1^*(\theta_a^L-\theta_a^R) 
-\f{1}{4}(e_a)^{2}_M \otimes\Phi_2^*(\theta_a^L-\theta_a^R).
$$
Next, using $\iota((e_a)^{i}_M) \omega =
\frac{1}{2} \Phi_{i}^* (\theta^L_a + \theta^R_a) \ , i= 1,2 \ ,$ we find 
that
$$
\frac{1}{2} \big((e_a)^{1}_M \wedge (e_a)^{2}_M\big)^\sharp\circ \om^\flat
=  \frac{1}{4} (e_a)^{1}_M \otimes \Phi^*_2 (\theta^L_a + \theta_a^R) 
- \frac{1}{4} (e_a)^{2}_M \otimes \Phi^*_1 (\theta^L_a + \theta_a^R).
$$
 From 
$P ^\sharp ( \Phi^*_{i} \theta^L_a) = \frac{1}{2} (1+
\Ad_{\Phi_i})_{ab} (e_b)^{i}_M \ , i =1,2 \ ,$ we obtain 
$$
\frac{1}{2} P^\sharp\circ \big( \Phi^*_1 \theta^L_a \wedge \Phi^*_2
\theta^R_a\big)^\flat= {\f {1}{4}}
(1+ \Ad_{\Phi_2^{-1}})_{ab} 
(e_a)^{2}_M \otimes \Phi^*_1 \theta^L_b
-{\f {1}{4}}(1+\Ad_{\Phi_1})_{ab} (e_a)^{1}_M \otimes 
\Phi^*_2 \theta^R_b \ .
$$
 Finally, 
\beq
{\f{1}{4}}((e_a)^{1}_M \wedge (e_a)^{2}_M)^\sharp 
& \circ & (\Phi^*_1 \theta^L_b  \wedge \Phi^*_2
\theta^R_b)^\flat \\
& = & 
{\f {1}{4}(e_a)^{1}_M \otimes (1- \Ad_{\Phi_2^{-1}})_{ab} \Phi^*_1 \theta^L_b +
{\f {1}{4}}(e_a)^{2}_M \otimes (\Ad_{\Phi_1}-1})_{ab} \Phi^*_2 \theta^R_b.
\eeq
Putting everything together, we get 
\beq 
P_{fus}^\sharp \circ \omega_{fus}^\flat &=& \Id_{TM} - \frac{1}{4}
(e_a)_M \otimes (\Phi^*_2 \theta^L_a + (\Ad_{\Phi_2^{-1}})_{ab}
\Phi^*_1 \theta^L_b - \Phi^*_1 \theta^R_a - (\Ad_{\Phi_1})_{ab}
\Phi_2^* \theta^R_b)\\
&=&
\Id_{TM} - \frac{1}{4}(e_a)_M
 \otimes
(\Phi_1 \Phi_2)^* (\theta^L_a -\theta_a^R) \, ,
\eeq
as required.
\end{proof}

\begin{example} \label{ex:moduli}
Consider
$$
M=\underbrace{{\bf D}(G)\fus\ldots 
\fus {\bf D}(G)}_{h}\fus \,  \Co_1 \fus \ldots \fus \Co_r,
$$
where $\Co_1, \dots, \Co_r$ are conjugacy classes in $G$.
This space can be viewed either as a quasi-Hamiltonian
space (see Section 9 of \cite{al:mom}), or as a quasi-Poisson manifold
(see Section 6). Proposition \ref{prop:fusioneq} together with 
Examples \ref{ex:Pconj} and
\ref{ex:Pdouble} show that these two structures agree in the sense of 
Theorem \ref{th:equivalence}. According to Theorem 9.3 of \cite{al:mom},
the reductions of $M$ are isomorphic to the moduli spaces of flat connections
with the Atiyah-Bott symplectic form. 
Proposition \ref{prop:redeq}
implies that the quasi-Poisson reduction yields the Poisson bivector inverse
to the canonical symplectic form.  
\end{example}

\begin{appendix} \label{app:KKS}
\section{The formal Poisson structure of $L\g^*$}
In this appendix, we show that the quasi-Poisson structure $P_{G}$  
on $G$ defined by Equation \eqref{eq:Gbivector}
can be viewed as a quotient of a formal  Poisson
structure on $L\g^*$, the dual of the loop algebra $L\g$ of $\g$. 

Let $LG$ be the loop group of $G$, let $\widehat{LG}$
be its central extension, and let $\widehat{L\g}^*$ be 
the dual of the Lie algebra $\widehat{L\g}$ of $\widehat{LG}$.
We let $L\g^*$ be the hyperplane at level 1 in $\widehat{L\g}^*$, 
equipped with the affine $LG$-action 
and the formal linear Poisson  bivector  
obtained by restriction from those on $\widehat{L\g}^*$. 

The action of the based loop group $\Omega G \subset LG$ on $L\g^*$ is free,
with quotient $L\g^*/\Omega G=G$. 
We shall show that,
under the quotient map, the formal Poisson
structure on $L\g^*$
projects to the bivector \eqref{eq:Gbivector} on $G$.

In the finite-dimensional case,  
the restriction of the linear Poisson structure $P_{0,\g^*}$ of $\g^{*}$
to the Lie algebra $\t$ of a 
maximal torus $T\subset G$ can be written
in terms of the corresponding root space 
decomposition. Let $\mf{R}$ be the 
set of roots of the complexified Lie algebra $\g^{\C}$, 
and for any $\alpha\in\mf{R}$, let 
$e_\alpha$ 
be a root vector of unit length, such that 
$e_{-\alpha}$ is the complex conjugate 
of $e_\alpha$ and $\langle e_\alpha, e_{-\alpha} \r =1$,
where $\l~,~\r$ is the Killing form. 
Then, if $\mu \in\t$ with $\l\alpha,\mu \r\not=0$ for all roots $\alpha$,
 the value of $P_{0,\g^*}$ at $\mu$ is
$$ P_{0,\g^*}=
\sum_{\alpha\in\mf{R}_+}
\f{1}{2\pi i\l\alpha,\mu \r }
(e_{\alpha})_{\g^*}\wedge (e_{-\alpha})_{\g^*} \ .
$$

In the case of loop algebras, for any $\xi\in\g^\C$ and $k\in\Z$, we 
denote by $\xi[k]\in L\g^\C$ the loop defined by
$$ \xi[k](e^{2\pi i s})=e^{2\pi i k s}\xi.$$
Then the root vectors for $\widehat{L\g}$ are 
all $e_\alpha[k]$, with corresponding affine roots $(\alpha,k)$, 
together with vectors $h_j[k]$ for $k\not=0$, with affine 
roots $(0,k)$, where $(h_j)$ is an orthonormal basis for $\t$. 
 From the identification of $L\g^*$ with the hyperplane at level $1$ in 
$\wh{L\g}^*$, we obtain the following formal expression 
for the Poisson structure on $L\g^*$ at a constant loop
$\mu\in\t\subset L\g^*$ 
with $\l\alpha,\mu\r\not\in\Z$, for all roots $\alpha$,  
\beq
P_{0,L\g^*} & = &
\sum_{\alpha\in\mf{R}_+}\sum_{k\in\Z}
\f{1}{2\pi i(\l\alpha,\mu \r+k)}(e_{\alpha}[k])_{L\g^*}\wedge 
(e_{-\alpha}[k])_{L\g^*} 
\\
& +&
 \sum_{j}
\sum_{k>0}\f{1}{2\pi i k}(h_j[k])_{L\g^*}\wedge (h_j[-k])_{L\g^*} \ .
\eeq

In a finite-dimensional manifold, any Poisson structure which is 
invariant under a free group action reduces to a Poisson 
structure on the quotient. Formally, we can carry out this calculation 
for the free action of $\Omega G$ on $L\g^*$.
The quotient map takes a constant loop
$\mu \in\t \subset L\g^*$ to the element $\exp\mu \in T\subset G$, 
and the corresponding tangent map takes
the value of $(\xi[k])_{L\g^*}$ at $\mu$ to the value of $\xi_G$ at 
${\exp\mu}$. Applying the tangent of the quotient map
to the Poisson bivector of 
$L\g^*$, we obtain 
$$ \sum_{\alpha\in\mf{R}_+} \sum_{k \in \Z}
\f{1}{2\pi i(\l\alpha,\mu\r+k)}(e_{\alpha})_G\wedge (e_{-\alpha})_G \, .
$$
Using the formula 
$$ \lim_{N\to \infty} \sum_{|k|\le N}\f{1}{2\pi (x+k)}
=\f{1}{2}\cot (\pi x)
$$
for $x\not\in\Z$, we obtain  
$$ \f{1}{2i}\sum_{\alpha\in\mf{R}_+}
\cot(\pi\l\alpha,\mu\r) (e_\alpha)_G\wedge (e_{-\alpha})_G \ ,
$$ 
which coincides with the bivector 
$P_G$ defined by Equation \eqref{eq:Gbivector} (see
Equation \eqref{eq:equationP}). 
An alternative procedure
for projecting the bivector of
$L\g^*$ to $G$ was considered in \cite{se:mo}.
\section{A generalized dynamical $r$-matrix}
\label{sec:lemma}
Let $r^{\g/\h}\in C^\infty(U,\wedge^2 \h^\perp)$ 
denote the $r$-matrix defined by Equation \eqref{eq:rmatrix}.
We shall now prove Lemma \ref{lem:YB}.  
Let $T\subset H$ be a maximal torus, and on $T\cap U$ define 
$$ r^{\g/\t}(h)=-\hh \Big(\f{\Ad_h+1}{\Ad_h-1}\Big|\t^\perp
\Big)_{ab}e_a\wedge e_b$$
and
$$  r^{\h/\t}(h)=- \hh \Big(\f{\Ad_h+1}{\Ad_h-1}\Big|\h\cap 
\t^\perp\Big)_{ab}e_a\wedge e_b \ .$$
Then $r^{\g/\h}|_{T\cap U}=r^{\g/\t}-r^{\h/\t}$. Our starting point will be 
the classical dynamical Yang-Baxter equations satisfied by 
$r^{\g/\t}$ and $r^{\h/\t}$ (see \cite{al:no}, Lemma A.5):
$$ 
\on{Cycl}_{abc}\big(\hh(e_a^L+e_a^R)_\t r_{bc}^{\g/\t}
+ r_{ak}^{\g/\t}f_{kbl}r_{lc}^{\g/\t} \big)=\f{1}{4}f_{abc}
$$  
and 
$$ 
\on{Cycl}_{abc}\big(\hh (e_a^L+e_a^R)_\t r_{bc}^{\h/\t}
+ r_{ak}^{\h/\t}f^\h_{kbl}r_{lc}^{\h/\t} \Big)=\f{1}{4}f^\h_{abc}.
$$ 
Using the fact that $r_{ak}^{\h/\t}f^\h_{kbl}r_{lc}^{\h/\t}=
r_{ak}^{\h/\t}f_{kbl}r_{lc}^{\h/\t}$, upon subtracting 
the second equation from the first, we obtain  
\begin{equation}
\label{eq:A}
\on{Cycl}_{abc}\big(\hh(e_a^L+e_a^R)_\t r_{bc}^{\g/\h}
+r_{ak}^{\g/\h}f_{kbl}r_{lc}^{\g/\h}
+2 r_{ak}^{\h/\t}f_{kbl}r_{lc}^{\g/\h} \big)=
\f{1}{4}(f_{abc}-f^\h_{abc}).
\end{equation}
To prove the Lemma, we need to evaluate
$\on{Cycl}_{abc}\big(\hh(e_a^L+e_a^R)_{\h/\t}
r_{bc}^{\g/\h}\big)$, where the subscript $\h/\t$ denotes the
projection of the vector field along $T\cap U$ onto the normal bundle
of $T\cap U$ in $U$ which is isomorphic to $(T\cap U)\times
(\t^\perp\cap \h)$. This projection can be expressed in terms of 
$r$ itself and the vector fields generating the $G$-action as
$$ \hh(e_a^L+e_a^R)_{\h/\t}=r_{ak}^{\h/\t} (e_k)_G.$$
By the $H$-invariance of $r^{\g/\h}$, this relation shows that
\begin{equation}\label{eq:B}
\on{Cycl}_{abc}\Big(\hh (e_a^L+e_a^R)_{\h/\t}r_{bc}^{\g/\h}
-2 r_{ak}^{\h/\t} f_{kbl} r_{lc}^{\g/\h}\Big)=0.
\end{equation}
Adding Equation \eqref{eq:A} to \eqref{eq:B}, we obtain Equation
\eqref{eq:EV2}, proving the Lemma.  
\end{appendix}

\bibliographystyle{amsplain}
%\bibliography{../../Biblio/ref}
%\bibliography{ref}

\providecommand{\bysame}{\leavevmode\hbox to3em{\hrulefill}\thinspace}

\end{document}